\documentclass[11pt,a4paper]{article}% For Computer Modern Fonts, or,
\usepackage{latexsym,amssymb, amsthm}
\usepackage[centertags]{amsmath}
\usepackage{epsfig,pifont}
\usepackage[inactive]{srcltx}  %SRC Specials for DVI search

\usepackage{graphics,helvet,times,mathptm,epic,eepic}

\usepackage{color}

\usepackage{newlfont}

\usepackage{amsfonts,bbm}

\newenvironment{keywords}{ \noindent {\small\bf Key Words}:}{ }

\def\bd{\begin{description}}
\def\ed{\end{description}}

\def\beq{\begin{equation}}
\def\eeq{\end{equation}}
\def\bea{\begin{eqnarray}}
\def\eea{\end{eqnarray}}
\def\beas{\begin{eqnarray*}}
\def\eeas{\end{eqnarray*}}

\newtheorem{theorem}{Theorem}[section]

\theoremstyle{remark}
\newtheorem{example}{Example}[section]

\begin{document}
%\begin{article}

\title{\textbf{\textsc{A new applied approach for  executing computations
 with infinite and infinitesimal quantities}}}

%% Author name
\newcommand{\nms}{\normalsize}
\author{  {   \bf Yaroslav D. Sergeyev\footnote{Yaroslav D. Sergeyev, Ph.D., D.Sc., holds a Full Professorship
 reserved for distinguished scientists at the University of Calabria, Rende, Italy.
 He is also Full Professor (a part-time contract) at the N.I.~Lobatchevsky State University,
  Nizhni Novgorod, Russia and Affiliated Researcher at the Institute of High Performance
  Computing and Networking of the National Research Council of Italy.   }
    }\\ \\ [-2pt]
      \nms Dipartimento di Elettronica, Informatica e Sistemistica,\\[-4pt]
       \nms   Universit\`a della Calabria,\\[-4pt]
       \nms 87030 Rende (CS)  -- Italy\\ \\[-4pt]
       \nms http://wwwinfo.deis.unical.it/$\sim$yaro\\[-4pt]
         \nms {\tt  yaro@si.deis.unical.it }
}

\date{}

\maketitle

\begin{abstract}
 A new computational methodology for executing calculations
 with infinite and infinitesimal quantities is described in
this paper. It is based on the principle `The part is less than
the whole' introduced by Ancient Greeks and applied to all numbers
(finite, infinite, and infinitesimal) and to all sets and
processes (finite and infinite). It is shown that it becomes
possible to write down finite, infinite, and infinitesimal numbers
by a finite number of symbols as particular cases of a unique
framework. The new methodology has allowed us to introduce the
Infinity Computer working with such numbers (its simulator has
already been realized).   Examples dealing with divergent series,
infinite sets, and limits   are given.
 \end{abstract}

\begin{keywords}
Infinite and infinitesimal numbers, infinite unite of measure,
numeral
    systems,  infinite sets, divergent series.
 \end{keywords}

\section{Introduction}
\label{s0}

Throughout the whole history of humanity many brilliant thinkers
studied problems related to the idea of infinity (see
\cite{Cantor,Cohen,Conway,Godel,Hardy,Hilbert,Leibniz,Newton,Robinson}
and references given therein). To emphasize importance of the
subject it is sufficient to mention that the Continuum Hypothesis
related to infinity has been included by David Hilbert as the
Problem Number One in his famous list of 23 unsolved mathematical
problems (see \cite{Hilbert}) that have influenced strongly
development of   Mathematics in the XX-th  century.

There exist different ways to generalize traditional arithmetic
for finite numbers to the case of infinite and infinitesimal
numbers (see \cite{Benci,Cantor,Conway,Robinson}  and references
given therein). However, arithmetics developed for infinite
numbers are quite different with respect to the finite arithmetic
we are used to deal with. Moreover, very often they leave
undetermined many operations where infinite numbers take part (for
example, $\infty-\infty$, $\frac{\infty}{\infty}$, sum of
infinitely many items, etc.) or use representation of infinite
numbers based on infinite sequences of finite numbers. In spite of
these crucial difficulties and due to enormous importance of the
concept of infinity in science, people try to introduce infinity
in their work with computers. We can mention the IEEE Standard for
Binary Floating-Point Arithmetic containing representations for
$+\infty$ and $-\infty$ and incorporation of these notions in the
interval analysis implementations.

%, and COSY system based on ideas of non-Archimedean fields (see
%\cite{Berz}).

The point of view on infinity accepted nowadays takes its origins
from the famous ideas of Georg Cantor (see \cite{Cantor}) who has
shown that there exist infinite sets having different number of
elements. However, it is well known that Cantor's approach leads
to some situations that often are called by non mathematicians
`paradoxes'. The most famous and simple of them is, probably,
Hilbert's paradox of the Grand Hotel. In a normal hotel having a
finite number   of rooms no more new guests can be accommodated if
it is full. Hilbert's Grand Hotel has an infinite number of rooms
(of course, the number of rooms is countable, because the rooms in
the Hotel are numbered). Due to Cantor, if a new guest arrives at
the Hotel where every room is occupied, it is, nevertheless,
possible to find a room for him. To do so, it is necessary to move
the guest occupying room 1 to room 2, the guest occupying room 2
to room 3, etc. In such a way room 1 will be ready for the
newcomer and, in spite of our assumption that there are no
available rooms in the Hotel, we have found one.

This result is very difficult to be fully realized by anyone who
is not a mathematician since in our every day experience in the
world around us the part is always less than the whole and if a
hotel is complete there are no places in it. In order to
understand how it is possible to tackle the problem of infinity in
such a way that Hilbert's Grand Hotel would be in accordance with
the principle `the part is less than the whole' let us consider a
study published in \textit{Science} by Peter Gordon (see
\cite{Gordon}) where he describes a primitive tribe living in
Amazonia - Pirah\~{a} - that uses a very simple numeral
system\footnote{ We remind that \textit{numeral}  is a symbol or
group of symbols that represents a \textit{number}. The difference
between numerals and numbers is the same as the difference between
words and the things they refer to. A \textit{number} is a concept
that a \textit{numeral} expresses. The same number can be
represented by different numerals. For example, the symbols `3',
`three', and `III' are different numerals, but they all represent
the same number.} for counting: one, two, many. For Pirah\~{a},
all quantities larger than two are just `many' and such operations
as 2+2 and 2+1 give the same result, i.e., `many'. Using their
weak numeral system Pirah\~{a} are not able to see, for instance,
numbers 3, 4, 5, and 6, to execute arithmetical operations with
them, and, in general, to say anything about these numbers because
in their language there are neither words nor concepts for that.
Moreover, the weakness of their numeral system leads to such
results as
\[
\mbox{`many'}+ 1= \mbox{`many'},   \hspace{1cm}    \mbox{`many'} +
2 = \mbox{`many'},
\]
which are very familiar to us  in the context of views on infinity
used in the traditional calculus
\[
\infty + 1= \infty,    \hspace{1cm}    \infty + 2 = \infty.
\]
This observation leads us to the following idea: \textit{Probably
our difficulty in working with infinity is not connected to the
nature of infinity but is a result of inadequate numeral systems
used to express numbers.}

In this paper, we describe a new methodology for treating infinite
and infinitesimal quantities (examples of its usage see in
\cite{Sergeyev,Sergeyev_patent,www,Philosophy,Poland,Mathesis,chaos})
having a strong numerical character. Its description is given in
Section~\ref{s1}. The new methodology allows us  to introduce in
Section~\ref{s2}  a new infinite unit of measure that is then used
as   the radix of a new positional numeral system.
Section~\ref{s3} shows that this system allows one to express
finite, infinite, and infinitesimal numbers in a unique framework
and to execute arithmetical operations with all of them.
Section~\ref{s4} discusses   some applications of the new
methodology. Section~\ref{s6} establishes relations to some of the
results of Georg Cantor. After all, Section~\ref{s7} concludes the
paper.

We close this Introduction by emphasizing   that the goal of the
paper is not to construct a complete theory of infinity and to
discuss such concepts as, for example, `set of all sets'. In
contrast, the problem of infinity is considered from the point of
view of applied Mathematics and theory and practice of
computations -- fields being among the main scientific interests
(see, e.g., monographs \cite{Sergeyev,Strongin_Sergeyev}) of the
author. A new viewpoint on infinity is introduced in the paper in
order to give possibilities to solve new and old (but with higher
precision) applied problems. Educational issues (see
\cite{Mockus_1,Mockus_2,Mathesis}) have also been taken into
account. In this connection, it is worthy to notice that a new
kind of computers -- the Infinity Computer -- able to execute
computations with infinite and infinitesimal numbers introduced in
this paper has been recently proposed and its software simulator
has already been implemented (see
\cite{Sergeyev_patent,www,Poland}).

%The Infinity Calculator using the Infinity Computer technology has
%been also implemented (see Fig.~\ref{Big_paper0}).
% \begin{figure}[ht]
%  \begin{center}
%    \epsfig{ figure = Big_paper0.eps, width = 3.2in, height = 2.6in,  silent = yes }
%    \caption{An example of multiplication at the Infinity Calculator.
%Arithmetical operations with finite, infinite, and infinitesimal
%numbers are discussed in Section~\ref{s3}.}
% \label{Big_paper0}
%  \end{center}
% \end{figure}

\section{A new computational methodology}
\label{s1}

The aim of this section is to introduce a new methodology that
would allow one to work with infinite and infinitesimal quantities
\textit{in the same way} as one works with finite numbers.
Evidently,  it becomes necessary to define what does it mean
\textit{in the same way}. Usually,  in modern Mathematics, when it
is necessary  to define a concept or an object, logicians   try to
introduce a number of axioms describing the object. However, this
way is fraught with danger because of the following reasons. First
of all, when we describe a mathematical object or concept we are
limited by the expressive capacity of the language we use to make
this description. A more rich language allows us to say more about
the object and a weaker language -- less (remind Pirah\~{a} that
are not able to say a word about number 4). Thus, development of
the mathematical (and not only mathematical) languages leads to a
continuous necessity of a transcription and specification of
axiomatic systems. Second, there is no any guarantee that the
chosen axiomatic system defines `sufficiently well' the required
concept and a continuous comparison with practice is required in
order to check the goodness of the accepted set of axioms.
However, there cannot be again any guarantee that the new version
will be the last and definitive one. Finally, the third limitation
  latent in axiomatic systems has been discovered by
G\"odel in his two famous incompleteness theorems (see
\cite{Godel_1931}).

In this paper, we introduce a different, significantly more
applied and less ambitious view on axiomatic systems related only
to utilitarian necessities to make calculations.   We start by
introducing three postulates that will fix our methodological
positions with respect to infinite and infinitesimal quantities
and Mathematics, in general. In contrast to the modern
mathematical fashion that tries to make all axiomatic systems more
and more precise (decreasing so degrees of freedom of the studied
part of Mathematics), we just define a set of general rules
describing how practical computations should be executed leaving
so as much space as possible for further, dictated by practice,
changes and developments of the introduced mathematical language.
Speaking metaphorically, we prefer to make a hammer and to use it
instead of describing what is a hammer and how it works.

Usually, when mathematicians deal with infinite objects (sets or
processes) it is supposed (even by constructivists (see, for
example, \cite{Markov})) that human beings are able to execute
certain operations infinitely many times. For example, in a fixed
numeral system it is possible to write down a numeral with
\textit{any} number of digits. However, this supposition is an
abstraction (courageously declared by constructivists in
\cite{Markov}) because we live in a finite world and all human
beings and/or computers finish operations they have started. In
this paper,    this abstraction is not used  and the following
postulate is adopted.

\textbf{Postulate 1.} \textit{We postulate  existence of infinite
and infinitesimal objects but accept that human beings and
machines are able to execute only a finite number of operations.}

Thus, we accept that we shall never be able to give a complete
description of infinite processes and sets due to our finite
capabilities.  Particularly, this means that we accept that we are
able to write down only a finite number of symbols to express
numbers.

The second postulate that will be adopted is due to the following
consideration. In natural sciences, researchers use tools to
describe the object of their study and the used instruments
influence results of observations. When physicists see a black dot
in their microscope they cannot say: the object of observation
\textit{is} the black dot. They are obliged to say: the lens used
in the microscope allows us to see the black dot and it is not
possible to say anything more about the nature of the object of
observation until we shall not change the instrument -- the lens
or the microscope itself -- by a more precise one.

Due to Postulate 1, the same happens in Mathematics when studying
natural phenomena, numbers, and objects that can be constructed by
using numbers. Numeral systems used to express numbers are among
the instruments of observations used by mathematicians. Usage of
powerful numeral systems gives the possibility to obtain more
precise results in mathematics in the same way as usage of a good
microscope gives the possibility to obtain more precise results in
Physics. However, the capabilities of the tools will be always
limited due to Postulate 1. Thus, following natural sciences, we
accept the second postulate.

\textbf{Postulate 2.} \textit{We shall not   tell \textbf{what
are} the mathematical objects we deal with; we just shall
construct more powerful tools that will allow us to improve our
capacities to observe and to describe properties of mathematical
objects.}

Particularly, this means that from this applied point of view,
axiomatic systems do not define mathematical objects but just
determine formal rules for operating with certain numerals
reflecting some properties of the studied mathematical objects.
For example,  axioms for real numbers are considered together with
a particular numeral system $\mathcal{S}$ used to write down
numerals and are viewed as practical rules (associative and
commutative properties of multiplication and addition,
distributive property of multiplication over addition, etc.)
describing operations with the numerals. The completeness property
is interpreted as a possibility to extend $\mathcal{S}$ with
additional symbols (e.g., $e$, $\pi$, $\sqrt{2}$, etc.) taking
care of the fact that the results of computations with these
symbols  agree with the facts observed in practice. As a rule, the
assertions regarding numbers that cannot be expressed in a numeral
system are avoided (e.g., it is not supposed that real numbers
form a field).

After all, we want to treat   infinite and infinitesimal numbers
in the same manner as we are used to deal with finite ones, i.e.,
by applying the philosophical principle of Ancient Greeks `The
part is less than the whole'. This principle, in our opinion, very
well reflects organization of the world around us but is not
incorporated in many traditional infinity theories where it is
true only for finite numbers.

\textbf{Postulate 3.} \textit{We adopt the principle `The part is
less than the whole' to all numbers (finite, infinite, and
infinitesimal) and to all sets and processes (finite and
infinite).}

Due to this declared applied statement, such concepts as
bijection, numerable and continuum sets, cardinal and ordinal
numbers cannot be used in this paper because they belong to
theories working with different assumptions. However, the approach
proposed here does not contradict Cantor. In contrast, it evolves
his deep ideas regarding existence of different infinite numbers
in a more applied way.

It is important to notice that the  adopted Postulates impose also
the style of  exposition of results in the paper: we first
introduce   new mathematical instruments, then show how to use
them in several areas of Mathematics, introducing each item as
soon as it becomes indispensable for the problem under
consideration.

Let us introduce now the main methodological idea of the paper  by
studying a situation  arising in practice and related to the
necessity to operate with extremely large quantities (see
\cite{Sergeyev} for a detailed discussion). Imagine that we are in
a granary and the owner asks us to count how much grain he has
inside it.  Of course, nobody counts the grain seed by seed
because the number of seeds is enormous.

To overcome this difficulty, people take sacks, fill them in with
seeds, and count the number of sacks. It is important that nobody
counts the number of seeds in a sack. If the granary is huge and
it becomes difficult to count the sacks, then trucks or even big
train waggons are used. Of course, we suppose that    all sacks
contain the same number of seeds, all trucks -- the same number of
sacks, and all waggons -- the same number of trucks. At the end of
the counting we obtain a result in the following form: the granary
contains 14 waggons, 54 trucks, 18 sacks, and 47 seeds of grain.
Note, that if we add, for example, one seed to the granary, we can
count it and see that the granary has more grain. If we take out
one waggon, we again are able to say how much grain has been
subtracted.

Thus, in our example it is  necessary to count large quantities.
They are finite but it is impossible to count them directly by
using an elementary unit  of measure, $u_0$, (seeds in our
example) because the quantities expressed in these units would be
too large. Therefore, people   are forced to behave as if the
quantities were infinite.

To solve the problem of `infinite' quantities, new units of
measure, $u_1,u_2,$ and $u_3,$ are introduced (units $u_1$ --
sacks, $u_2$ -- trucks, and $u_3$ -- waggons). The new units have
the following important peculiarity:  all the units $u_{i+1}$
contain a certain number $K_i$\label{p:1} of units $u_{i}$ but
this number, $K_i$, is unknown. Naturally, it is supposed that
$K_i$ is the same for all instances of the units  $u_{i+1}$. Thus,
numbers that it was impossible to express using only the initial
unit  of measure are perfectly expressible in the new units we
have introduced in spite of the fact that the numbers  $K_i$ are
unknown.

This key idea of counting by introduction of   new units of
measure will be used in the paper to deal with infinite quantities
together  with the idea of separate  count of units with different
exponents used in traditional positional numeral systems.

\section{The infinite unit of measure}
\label{s2}

The infinite unit of measure is expressed by the numeral
\ding{172} called   \textit{grossone} and is introduced as the
number of elements of  the set, $\mathbb{N}$, of natural numbers.
Remind that the usage of a numeral indicating totality of the
elements we deal with is not new in Mathematics. It is sufficient
to mention the theory of probability (axioms of Kolmogorov) where
events can be defined in two ways. First, as union of elementary
events; second, as a sample space, $\Omega$, of all possible
elementary events (or its parts $\Omega/2, \Omega/3,$ etc.) from
which some elementary events have been excluded (or added in case
of parts of $\Omega$). Naturally, the latter way to define events
becomes particularly useful when the sample space consists of
infinitely many elementary events.

Grossone  is introduced by describing its properties  (similarly,
in order to pass from natural to integer numbers a new element --
zero -- is introduced by describing its properties) postulated by
the \textit{Infinite Unit Axiom} (IUA) consisting of three parts:
Infinity, Identity, and Divisibility. This axiom is added to
axioms for real numbers (remind that we consider axioms in sense
of Postulate~2). Thus, it is postulated that associative and
commutative properties of multiplication and addition,
distributive property of multiplication over addition, existence
of   inverse  elements with respect to addition and multiplication
hold for grossone as for finite numbers\footnote{It is important
to emphasize that we speak about axioms of real numbers in sense
of Postulate~2, i.e., axioms define formal rules of operations
with numerals in a given numeral system. Therefore, if we want to
have a numeral system including grossone, we should fix also a
numeral system to express finite numbers. In order to concentrate
our attention on properties of grossone, this point will be
investigated later. }. Let us introduce the axiom and then give
comments on it.

\textit{Infinity.}
  Any finite  natural number $n$  is less than grossone, i.e.,  $n
<~\mbox{\ding{172}}$.

\textit{Identity.}
 The following
relations  link \ding{172} to identity elements 0 and 1
 \beq
 0 \cdot \mbox{\ding{172}} =
\mbox{\ding{172}} \cdot 0 = 0, \hspace{3mm}
\mbox{\ding{172}}-\mbox{\ding{172}}= 0,\hspace{3mm}
\frac{\mbox{\ding{172}}}{\mbox{\ding{172}}}=1, \hspace{3mm}
\mbox{\ding{172}}^0=1, \hspace{3mm}
1^{\mbox{\tiny{\ding{172}}}}=1, \hspace{3mm}
0^{\mbox{\tiny{\ding{172}}}}=0.
 \label{3.2.1}
       \eeq

\textit{Divisibility.}
For any finite natural number  $n$   sets $\mathbb{N}_{k,n}, 1 \le
k \le n,$ being the $n$th parts of the set, $\mathbb{N}$, of
natural numbers have the same number of elements indicated by the
numeral $\frac{\mbox{\ding{172}}}{n}$ where
 \beq
 \mathbb{N}_{k,n} = \{k,
k+n, k+2n, k+3n, \ldots \}, \hspace{5mm} 1 \le k \le n,
\hspace{5mm} \bigcup_{k=1}^{n}\mathbb{N}_{k,n}=\mathbb{N}.
 \label{3.3}
       \eeq

The first part of the introduced axiom -- Infinity -- is quite
clear. In fact, we want to describe an infinite number, thus, it
should be larger than any finite number. The second part of the
axiom -- Identity -- tells us that \ding{172}  behaves itself with
identity elements 0 and 1 as all other numbers. In reality, we
could even omit this part of the axiom because, due to
Postulate~3, all numbers should be treated in the same way and,
therefore, at the moment we have told that grossone is a number,
we have fixed usual properties of numbers, i.e., the properties
described in Identity, associative and commutative properties of
multiplication and addition, distributive property of
multiplication over addition, existence of inverse elements with
respect to addition and multiplication. The third part of the
axiom -- Divisibility -- is the most interesting, it is based on
Postulate~3. Let us first illustrate it by an example.

\begin{example}
\label{e0} If we take $n = 1$, then $\mathbb{N}_{1,1} =
\mathbb{N}$ and Divisibility tells that the set, $\mathbb{N}$, of
natural numbers has \ding{172} elements. If $n = 2$, we have two
sets $\mathbb{N}_{1,2}$ and $\mathbb{N}_{2,2}$
 \beq
\begin{array}{ccccccccccc}
     \mathbb{N}_{1,2} = &\hspace{-2mm} \{1, &   & 3, &  & 5, &   & 7, & \ldots  &\} \\
      &   &   &   &   &   &   &    &  \\
     \mathbb{N}_{2,2} = &  \hspace{-5mm}  \{  & 2, &  & 4, &  & 6, &   &  \ldots &\}  \\
    \end{array}
\label{3.3.1}
       \eeq
and they have $\frac{\mbox{\ding{172}}}{2}$  elements each. If $n
= 3$, then we have three sets
 \beq
\begin{array}{cccccccccc}
     \mathbb{N}_{1,3} = &\hspace{-2mm}\{1, &   &   & 4, &   &   & 7, &  \ldots   &\} \\
         &   &   &   &   &   &   &    &  \\
      \mathbb{N}_{2,3} = &    \hspace{-5mm}        \{  & 2, &   &   & 5, &   &   &  \ldots &\}     \\
                     &   &   &   &   &   &   &    &  \\
       \mathbb{N}_{3,3} = &     \hspace{-5mm}      \{  &   & 3, &   &   & 6, &   &  \ldots  &\}  \\
    \end{array}\vspace{ 2mm}
\label{3.3.2}
       \eeq
and they have $\frac{\mbox{\ding{172}}}{3}$  elements each.
\hfill$\Box$ \end{example}

It is important to emphasize that to introduce
$\frac{\mbox{\ding{172}}}{n}$ we do not try to count elements $k,
k+n, k+2n, k+3n, \ldots$ one by one in (\ref{3.3}). In fact, we
cannot do this due to Postulate~1.  By using Postulate~3, we
construct the sets $\mathbb{N}_{k,n}, 1 \le k \le n,$ by
separating the whole, i.e., the set $\mathbb{N}$, in $n$ parts
(this separation is highlighted visually in formulae (\ref{3.3.1})
and (\ref{3.3.2})). Again due to Postulate~3, we affirm that the
number of elements of the $n$th part of the set, i.e.,
$\frac{\mbox{\ding{172}}}{n}$, is $n$ times less than the number
of elements of the whole set, i.e., than \ding{172}. In terms of
our granary example  \ding{172} can be interpreted as the number
of seeds in the sack. Then, if the sack contains \ding{172} seeds,
its $n$th part contains $n$ times less quantity, i.e.,
$\frac{\mbox{\ding{172}}}{n}$ seeds. Note  that, since the numbers
$\frac{\mbox{\ding{172}}}{n}$ have been introduced as numbers  of
elements of   sets $\mathbb{N}_{k,n}$, they are integer.

The new unit of measure allows us to calculate easily the number
of elements of sets being union, intersection, difference, or
product of other sets of the type $\mathbb{N}_{k,n}$.  Due to our
accepted methodology, we do it  in the same way  as these
measurements are executed for finite sets. Let us consider two
simple examples (a general rule for determining the number of
elements of infinite sets having a more complex structure will be
given in Section~\ref{s4}) showing how grossone can be used for
this purpose.
\begin{example}
\label{e1}  Let us determine the number of elements of the set
$A_{k,n} = \mathbb{N}_{k,n} \backslash  \{a\},$ $a \in
\mathbb{N}_{k,n}, n \ge 1$. Due to the IUA, the set
$\mathbb{N}_{k,n}$ has $\frac{\mbox{\ding{172}}}{n}$ elements. The
set $A_{k,n}$ has been constructed by excluding one element from
$N_{k,n}$. Thus, the set $A_{k,n}$ has
$\frac{\mbox{\ding{172}}}{n}-1$ elements. The granary
interpretation can be also given for the number
$\frac{\mbox{\ding{172}}}{n}-1$: the number of seeds in the $n$th
part of the sack minus one seed. For $n=1$ we have
$\mbox{\ding{172}}-1$ interpreted as the number of seeds in the
sack minus one seed. \hfill$\Box$
 \end{example}
\begin{example}
\label{e2} Let us consider the following two sets
\[
B_1 = \{ 4, 9, 14, 19, 24, 29, 34, 39, 44, 49, 54, 59, 64, 69, 74,
79,\ldots\},
\]
\[
B_2 = \{ 3, 14, 25, 36, 47, 58, 69, 80, 91, 102, 113, 124, 135,
\ldots\}
\]
and determine the number of elements   in the set $B  = (B_1 \cap
B_2 ) \cup \{  3,4,5, 69  \}$. It follows immediately from the IUA
that $B_1 = \mathbb{N}_{4,5}, B_2 = \mathbb{N}_{3,11}$. Their
intersection
\[
B_1 \cap B_2 = \mathbb{N}_{4,5} \cap \mathbb{N}_{3,11} = \{ 14,
69, 124, \ldots\} = \mathbb{N}_{14,55}
\]
and, therefore, due to the IUA, it has
$\frac{\mbox{\ding{172}}}{55}$ elements. Finally, since 69 belongs
to the set  $\mathbb{N}_{14,55}$ and 3, 4, and 5 do not belong to
it, the set $B$ has $\frac{\mbox{\ding{172}}}{55}+3$ elements. The
granary interpretation: this is the number of seeds in the $55$th
part of the sack plus three seeds. \hfill$\Box$
\end{example}

One of the important differences of the new approach with respect
to the non-standard analysis consists of the fact that
$\mbox{\ding{172}} \in \mathbb{N}$ because grossone has been
introduced as the quantity of natural numbers  (similarly, the
number 5 being the number of elements of the set $\{1, 2, 3, 4, 5
\}$ is the largest element in this set). The new numeral
\ding{172} allows one to write down the set, $\mathbb{N}$, of
natural numbers in the form
 \beq
\mathbb{N} = \{ 1,2,3, \hspace{5mm} \ldots  \hspace{5mm}
 \mbox{\ding{172}}-3, \hspace{2mm} \mbox{\ding{172}}-2,
\hspace{2mm}\mbox{\ding{172}}-1, \hspace{2mm} \mbox{\ding{172}} \}
\label{4.1}
       \eeq
   where the numerals
 \beq
\ldots  \hspace{2mm} \mbox{\ding{172}}-3,
\hspace{2mm}\mbox{\ding{172}}-2, \hspace{2mm}\mbox{\ding{172}}-1,
\hspace{2mm} \mbox{\ding{172}} \label{4.2}
       \eeq
 indicate \textit{infinite} natural numbers.

It is important to emphasize that in the new approach the set
(\ref{4.1}) is the same set of natural numbers
 \beq
\mathbb{N} = \{ 1,2,3, \hspace{2mm} \ldots \hspace{2mm}  \}
\label{4.1_calcolo}
       \eeq
 we
are used to deal with and infinite numbers (\ref{4.2}) also take
part of $\mathbb{N}$. Both records, (\ref{4.1}) and
(\ref{4.1_calcolo}), are correct and do not contradict each other.
They just use two different numeral systems to express
$\mathbb{N}$. Traditional numeral systems   do not allow us to see
infinite natural numbers that we can observe now thanks
to~\ding{172}. Similarly,   Pirah\~{a}  are not able to see finite
natural numbers greater than~2. In spite of this fact, these
numbers (e.g., 3 and 4) belong to $\mathbb{N}$ and are visible if
one uses a more powerful numeral system. Thus, we have the same
object of observation -- the set $\mathbb{N}$ -- that can be
observed by different instruments -- numeral systems -- with
different accuracies
 (see Postulate~2).

Now the following obvious  question arises: Which natural numbers
can we express by using the new numeral \ding{172}? Suppose that
we have a numeral system, $\mathcal{S}$, for expressing finite
natural numbers and it allows us  to express $K_{\mathcal{S}}$
numbers (not necessary consecutive)  belonging to a set
$\mathcal{N}_{\mathcal{S}} \subset \mathbb{N}$. Note that due to
Postulate~1, $K_{\mathcal{S}}$ is finite. Then,   addition of
\ding{172} to this numeral system will allow us to express also
infinite natural numbers $\frac{i\mbox{\small{\ding{172}}}}{n}\pm
k \le \mbox{\ding{172}}$ where  $1 \le i \le n,\,\,\, k\in
\mathcal{N}_{\mathcal{S}},\,\,\, n \in \mathcal{N}_{\mathcal{S}}$
(note that since $\frac{\mbox{\small{\ding{172}}}}{n}$ are
integers, $\frac{i\mbox{\small{\ding{172}}}}{n}$ are integers
too). Thus, the more powerful system $\mathcal{S}$ is used to
express finite numbers, the more infinite numbers can be expressed
but their quantity  is always finite, again due to Postulate~1.
The new numeral system using grossone allows us to express more
numbers than traditional numeral systems thanks to the introduced
new numerals but, as it happens for all numeral systems, its
abilities to express numbers are   limited.

\begin{example} \label{e3}  Let us consider the  numeral system,
$\mathcal{P}$, of Pirah\~{a} able to express only numbers 1 and 2
(the only difference will be in the usage of  numerals `1' and `2'
instead of original numerals $I$ and $II$ used by Pirah\~{a}). If
we add to $\mathcal{P}$ the new numeral \ding{172}, we obtain a
new numeral system (we call it $\widehat{\mathcal{P}}$) allowing
us to express only ten numbers represented by the following
numerals
 \beq \underbrace{1,2}_{finite},
\hspace{5mm} \ldots \hspace{5mm}
\underbrace{\frac{\mbox{\small{\ding{172}}}}{2}-2,
\frac{\mbox{\small{\ding{172}}}}{2}-1,
\frac{\mbox{\small{\ding{172}}}}{2},
\frac{\mbox{\small{\ding{172}}}}{2}+1,
\frac{\mbox{\small{\ding{172}}}}{2}+2}_{infinite}, \hspace{5mm}
\ldots \hspace{5mm} \underbrace{\mbox{\ding{172}}-2,
\mbox{\ding{172}}-1, \mbox{\ding{172}}}_{infinite}.
 \label{4.2.1}
       \eeq
The first two numbers in (\ref{4.2.1}) are finite, the remaining
eight are infinite, and dots show  natural numbers that are not
expressible in $\widehat{\mathcal{P}}$. As a consequence,
$\widehat{\mathcal{P}}$ does not allow us to execute such
operation as $2+2$ or to add $2$ to
$\frac{\mbox{\small{\ding{172}}}}{2}+2$ because their results
cannot be expressed in it.  Of course, we do not say that results
of these operations are equal (as Pirah\~{a} do for operations
   $2+2$ and $2+1$). We just say that the results are not
expressible in $\widehat{\mathcal{P}}$ and it is necessary to take
another, more powerful numeral system if we want to execute these
operations. \hfill$\Box$ \end{example}

Note that   crucial limitations discussed in Example~\ref{e3} hold
for   sets, too. As a consequence, the numeral system
$\mathcal{P}$ allows us to define only the sets $\mathbb{N}_{1,2}$
and $\mathbb{N}_{2,2}$ among all possible sets of the form
$\mathbb{N}_{k,n}$ from (\ref{3.3}) because we have only two
finite numerals, `1' and `2', in $\mathcal{P}$. This numeral
system is too weak to define other sets of this type because
numbers greater than 2 required for these definition are not
expressible in $\mathcal{P}$. These limitations have a general
character and are related to all questions requiring a numerical
answer (i.e., an answer expressed only in numerals, without
variables). In order to obtain such an answer, it is necessary to
know at least one numeral system able to express numerals required
to write down this answer.

% For example, we are not able to answer to the following
% question: What is the number of \textit{all} the sets of the type
% $\mathbb{N}_{k,n}$? But  for \textit{any known} numeral system
% $\mathcal{S}$ we can give a clear numerical answer to the question
% being an applied version of the given above theoretical problem:
% What is the number of the sets of the type $\mathbb{N}_{k,n}$
% expressible in $\mathcal{S}$? If the fixed system $\mathcal{S}$ is
% able to express $P
% > 2$ finite natural numbers, then it allows us to define
% $[P(P+2)/2]$ sets $\mathbb{N}_{k,n}, 2 \le n \le P,$ where $[u]$
% is the integer part of $u$.

We are ready now to formulate the following important result being
a direct consequence of the accepted methodological postulates.
\begin{theorem}
\label{t1} The set $\mathbb{N}$ is not a monoid under addition.
\end{theorem}
\textit{Proof.}
 Due to Postulate~3, the operation $\mbox{\ding{172}}+1$ gives us as
the  result a number greater than \ding{172}. Thus, by definition
of grossone, $\mbox{\ding{172}}+1$ does not belong to $\mathbb{N}$
and, therefore, $\mathbb{N}$ is not closed under addition and is
not a monoid. \hfill$\Box$

This result also means that adding  the IUA to the axioms of
natural numbers    defines the set of \textit{extended natural
numbers} indicated as $\widehat{\mathbb{N}}$ and including
$\mathbb{N}$ as a proper subset
 \beq
  \widehat{\mathbb{N}} = \{
1,2, \ldots ,\mbox{\ding{172}}-1, \mbox{\ding{172}},
\mbox{\ding{172}}+1, \ldots , \mbox{\ding{172}}^2-1,
\mbox{\ding{172}}^2, \mbox{\ding{172}}^2+1, \ldots \}.
\label{4.2.2}
       \eeq
The extended natural numbers  greater than grossone are also
linked to sets of numbers and can   be interpreted in the terms of
grain.
\begin{example}
\label{e4} Let us determine the number of elements of the set
 \[
 C  =
\{
  (a_1, a_2, \ldots, a_m ) : a_i \in   \mathbb{N}, 1 \le i \le m
 \}.
 \]
The  elements of $C$ are $m$-tuples of natural numbers. It is
known from combinatorial calculus that if we have $m$ positions
and each of them can be filled in by one of $l$ symbols, the
number of the obtained $m$-tuples is equal to $l^m$. In our case,
since   $\mathbb{N}$  has grossone elements, $l =
\mbox{\ding{172}}$. Thus, the set $C$ has $\mbox{\ding{172}}^m$
elements. The granary interpretation: if we accept that the
numbers $K_i$ from page~\pageref{p:1} are such that
$K_i=\mbox{\ding{172}}, 1 \le i \le m-1,$ then
$\mbox{\ding{172}}^2$ can be viewed as the number of seeds in the
truck, $\mbox{\ding{172}}^3$ as the number of seeds in the train
waggon, etc. \hfill$\Box$
\end{example}

The set, $\widehat{\mathbb{Z}}$, of \textit{extended integer
numbers}\index{extended integer numbers} can be construct\-ed from
the set, $\mathbb{Z}$, of integer numbers by a complete analogy
and inverse elements with respect to addition are introduced
naturally. For example, $7\mbox{\ding{172}}$ has its inverse with
respect to addition equal to~$-7\mbox{\ding{172}}$.

It is important to notice that,  due to Postulates 1 and 2, the
new system of counting cannot give answers to \textit{all}
questions regarding infinite sets. What can we say, for instance,
about the number of elements of the sets $\widehat{\mathbb{N}}$
and $\widehat{\mathbb{Z}}$? The introduced numeral system based on
\ding{172} is too weak to give  answers to these questions. It is
necessary to introduce in a way a more powerful numeral system by
defining new numerals (for instance, \ding{173}, \ding{174}, etc).

We conclude this section by the following  remark. The IUA
introduces a new number -- the quantity of elements in the set of
natural numbers -- expressed by the new numeral \ding{172}.
However, other numerals and sets can be used to state the idea of
the axiom. For example, the numeral \ding{182}  can be introduced
as the number of elements of the set, $\mathbb{E}$, of even
numbers and can be taken as the base of a numeral system. In this
case, the IUA can be reformulated using the numeral \ding{182} and
numerals using it will be used to express infinite numbers. For
example, the number of elements of the set, $\mathbb{O}$, of odd
numbers will be expressed as $|\mathbb{O}|=|\mathbb{E}|=$
\ding{182} and $|\mathbb{N}|=2 \cdot$ \ding{182}. We emphasize
through this note that infinite  numbers (similarly to the finite
ones) can be expressed by various numerals and in different
numeral systems.

\section{Arithmetical operations in the new numeral system}
\label{s3}

We  have already started to write down simple infinite numbers and
to execute arithmetical operations with them without concentrating
our attention upon this question. Let us consider it
systematically.

\subsection{Positional numeral system with infinite radix}

Different numeral systems have been developed  to describe finite
numbers. In positional numeral systems, fractional numbers are
expressed by the record
 \beq
 (a_{n}a_{n-1} \ldots a_1 a_0 .
 a_{-1} a_{-2}  \ldots   a_{-(q-1)}     a_{-q})_b
 \label{3.10}
       \eeq
where numerals $a_i, -q \le i \le n,$ are called \textit{digits},
belong to the alphabet $\{ 0, 1, \ldots , b-1 \}$, and the dot is
used to separate  the fractional part from the integer one. Thus,
the numeral (\ref{3.10}) is equal to the   sum
 \beq
 a_{n} b^{n} + a_{n-1} b^{n-1} +  \ldots + a_1 b^1 +a_0 b^0+
 a_{-1} b^{-1} +   \ldots  + a_{-(q-1)}b^{-(q-1)} +  a_{-q}
 b^{-q}.
\label{3.11}
       \eeq
Record    (\ref{3.10})   uses numerals consisting of one symbol
each, i.e., digits $a_{i} \in \{ 0, 1,$ $ \ldots , b-1 \}$, to
express how many finite units of   the type $b^{i}$ belong to the
number (\ref{3.11}). Quantities of finite units $b^{i}$ are
counted separately for each exponent $i$ and all symbols in the
alphabet $\{ 0, 1, \ldots , b-1 \}$ express finite numbers.

To express infinite and infinitesimal numbers we shall use records
that are similar to (\ref{3.10}) and (\ref{3.11}) but have some
peculiarities. In order to construct a  number $C$   in the new
numeral positional system with base \ding{172}, we subdivide $C$
into groups corresponding to powers of \ding{172}:
 \beq
  C = c_{p_{m}}
\mbox{\ding{172}}^{p_{m}} +  \ldots + c_{p_{1}}
\mbox{\ding{172}}^{p_{1}} +c_{p_{0}} \mbox{\ding{172}}^{p_{0}} +
c_{p_{-1}} \mbox{\ding{172}}^{p_{-1}}   + \ldots   + c_{p_{-k}}
 \mbox{\ding{172}}^{p_{-k}}.
\label{3.12}
       \eeq
 Then, the record
 \beq
  C = c_{p_{m}}
\mbox{\ding{172}}^{p_{m}}    \ldots   c_{p_{1}}
\mbox{\ding{172}}^{p_{1}} c_{p_{0}} \mbox{\ding{172}}^{p_{0}}
c_{p_{-1}} \mbox{\ding{172}}^{p_{-1}}     \ldots c_{p_{-k}}
 \mbox{\ding{172}}^{p_{-k}}
 \label{3.13}
       \eeq
represents  the number $C$, where all numerals $c_i\neq0$, they
belong  to a traditional numeral system and are called
\textit{grossdigits}. They express finite positive or negative
numbers and show how many corresponding units
$\mbox{\ding{172}}^{p_{i}}$ should be added or subtracted in order
to form the number $C$. Grossdigits can be expressed by several
symbols using positional systems, the form $\frac{Q}{q}$ where $Q$
and $q$ are integer numbers, or in any other finite numeral
system.

Numbers $p_i$ in (\ref{3.13}) called \textit{grosspowers}  can be
finite, infinite, and infinitesimal (the introduction of
infinitesimal numbers will be given soon), they   are sorted in
the decreasing order
\[
p_{m} >  p_{m-1}  > \ldots    > p_{1} > p_0 > p_{-1}  > \ldots
p_{-(k-1)}  >   p_{-k}
 \]
with $ p_0=0$.

 In the traditional record
(\ref{3.10}), there exists a convention that a digit  $a_i$ shows
how many powers $b^i$ are present in the number and the radix $b$
is not written explicitly. In the record (\ref{3.13}), we write
$\mbox{\ding{172}}^{p_{i}}$ explicitly because in the new numeral
positional system  the number   $i$ in general is not equal to the
grosspower $p_{i}$. This gives possibility to write, for example,
such a number  as
$7.6\mbox{\ding{172}}^{244.5}\,34\mbox{\ding{172}}^{32}$ having
grospowers $p_2=244.5,p_{1}=32$ and grossdigits $c_{244.5}=7.6,
c_{32} = 34$ without indicating grossdigits equal to zero
corresponding to grosspowers less than 244.5 and greater than 32.
Note also that if a grossdigit $c_{p_i}=1$ then we often write
$\mbox{\ding{172}}^{p_i}$ instead of $1\mbox{\ding{172}}^{p_i}$.

\textit{Finite numbers} in this new numeral system are represented
by numerals having only one grosspower $ p_0=0$. In fact, if we
have a number $C$ such that $m=k=$~0 in representation
(\ref{3.13}), then due to (\ref{3.2.1}),   we have $C=c_0
\mbox{\ding{172}}^0=c_0$. Thus, the number $C$ in this case does
not contain grossone and is equal to the grossdigit $c_0$ being a
conventional finite number     expressed in a traditional finite
numeral system.

\textit{Infinitesimal numbers}   are represented by numerals $C$
having only negative finite or infinite grosspowers. The following
two numbers are examples of infinitesimals:
$3\mbox{\ding{172}}^{-3.2}$,
$37\mbox{\ding{172}}^{-2}11\mbox{\ding{172}}^{-15}$. The simplest
infinitesimal number is
$\mbox{\ding{172}}^{-1}=\frac{1}{\mbox{\ding{172}}}$ being the
inverse element with respect to multiplication for \ding{172}:
 \beq
\frac{1}{\mbox{\ding{172}}}\cdot\mbox{\ding{172}}=\mbox{\ding{172}}\cdot\frac{1}{\mbox{\ding{172}}}=1.
 \label{3.15.1}
       \eeq
Note that all infinitesimals are not equal to zero. Particularly,
$\frac{1}{\mbox{\ding{172}}}>0$ because it is a result of division
of two positive numbers.  It also has a clear granary
interpretation. Namely, if we have a sack   containing \ding{172}
seeds, then one sack divided by the number of seeds in it is equal
to one seed. Vice versa, one seed, i.e.,
$\frac{1}{\mbox{\ding{172}}}$, multiplied   by the number of seeds
in the sack, $\mbox{\ding{172}}$, gives one sack of seeds. Note
that the usage of infinitesimals as grosspowers can lead to more
complex constructions, particularly, again to infinitesimals, see
e.g., the number
$1\mbox{\ding{172}}^{\mbox{\tiny\ding{172}}^{-1}}(-1)\mbox{\ding{172}}^0$.

\textit{Infinite numbers}  in this   numeral system  are expressed
by numerals having at least one finite or infinite gross\-power
greater than zero. Thus, they have infinite parts and can also
have  a finite part and infinitesimal ones. If power
$\mbox{\ding{172}}^{0}$ is the lowest in a number then we often
write simply   grossdigit $c_0$ without $\mbox{\ding{172}}^{0}$,
for instance, we write $23\mbox{\ding{172}}^{14}5$ instead of
$23\mbox{\ding{172}}^{14}5\mbox{\ding{172}}^{0}$.
\begin{example}
\label{e5} The left-hand expression below shows how to write down
numbers in the new numeral system and the right-hand shows how the
value of the number is calculated:
\[
 15\mbox{\ding{172}}^{1.4\mbox{\tiny{\ding{172}}}}(-17.2045)\mbox{\ding{172}}^{3}7\mbox{\ding{172}}^{0}52.1\mbox{\ding{172}}^{-6}
=
15\mbox{\ding{172}}^{1.4\mbox{\tiny{\ding{172}}}}-17.2045\mbox{\ding{172}}^{3}+7\mbox{\ding{172}}^{0}+52.1\mbox{\ding{172}}^{-6}.
\]
The number above has one infinite part having the infinite
grosspower, one infinite part having the  finite grosspower, a
finite part, and an infinitesimal part. \hfill$\Box$
\end{example}

Finally, numbers having a finite and infinitesimal parts   can be
also expressed in the new numeral system, for instance, the number
$-3.5\mbox{\ding{172}}^{0}(-37)\mbox{\ding{172}}^{-2}11\mbox{\ding{172}}^{-15\mbox{\tiny{\ding{172}}}+2.3}$
has a finite  and two infinitesimal parts, the second  of them has
the infinite negative gross\-power equal to
$-15\mbox{\ding{172}}+2.3$.

\subsection{Arithmetical operations }

We start the description of arithmetical operations for the new
positional numeral system by the operation of \textit{addition}
(\textit{subtraction} is a direct consequence of addition and is
thus omitted) of two given infinite numbers $A$ and $B$, where
 \beq
A= \sum_{i=1}^{K} a_{k_{i}}\mbox{\ding{172}}^{k_{i}}, \hspace{1cm}
B= \sum_{j=1}^{M} b_{m_{j}}\mbox{\ding{172}}^{m_{j}}, \hspace{1cm}
C= \sum_{i=1}^{L} c_{l_{i}}\mbox{\ding{172}}^{l_{i}},
 \label{3.20}
       \eeq
and the result $C=A+B$ is constructed    by including in it all
items $a_{k_{i}}\mbox{\ding{172}}^{k_{i}}$ from $A$ such that
$k_{i} \neq m_{j},1 \le j \le M,$ and all items
$b_{m_{j}}\mbox{\ding{172}}^{m_{j}}$ from $B$ such that $m_{j}
\neq k_{i},1 \le i \le K$. If in $A$ and $B$ there are items such
that $k_{i}=m_{j}$, for some $i$ and $j$, then this grosspower
$k_{i}$ is included in $C$ with the grossdigit
$b_{k_{i}}+a_{k_{i}}$, i.e., as
$(b_{k_{i}}+a_{k_{i}})\mbox{\ding{172}}^{k_{i}}$.

\begin{example}
\label{e6}
We consider two infinite numbers $A$ and $B$, where
$$
A=16.5\mbox{\ding{172}}^{44.2}(-12)\mbox{\ding{172}}^{12}
17\mbox{\ding{172}}^{0}, \hspace{1cm} B=6.23\mbox{\ding{172}}^{3}
10.1\mbox{\ding{172}}^{0}15\mbox{\ding{172}}^{-4.1}.
$$
Their sum $C$ is calculated as follows:
\[
C=A+B=16.5\mbox{\ding{172}}^{44.2}+(-12)\mbox{\ding{172}}^{12}+
17\mbox{\ding{172}}^{0}+ 6.23\mbox{\ding{172}}^{3}+
10.1\mbox{\ding{172}}^{0} +15\mbox{\ding{172}}^{-4.1}=
\]
\[
16.5\mbox{\ding{172}}^{44.2}-
12\mbox{\ding{172}}^{12}+6.23\mbox{\ding{172}}^{3}+
27.1\mbox{\ding{172}}^{0}+15\mbox{\ding{172}}^{-4.1}=
 \]
\[
\begin{tabular}{cr}\hspace {28mm}$16.5\mbox{\ding{172}}^{44.2}
(-12)\mbox{\ding{172}}^{12}6.23\mbox{\ding{172}}^{3}
27.1\mbox{\ding{172}}^{0}15\mbox{\ding{172}}^{-4.1}.\hfill { }$ &
\hspace {2cm}
 $\Box$
\end{tabular}
\]
\end{example}

 The operation of \textit{multiplication}  of two numbers
 $A$ and $B$ in the form (\ref{3.20}) returns, as
the result, the infinite number $C$ constructed as follows:
 \beq
C= \sum_{j=1}^{M} C_{j}, \hspace{5mm} C_{j} =
b_{m_{j}}\mbox{\ding{172}}^{m_{j}}\cdot A =\sum_{i=1}^{K}
a_{k_{i}}b_{m_{j}}\mbox{\ding{172}}^{k_{i}+m_{j}}, \hspace{5mm}1
\le j \le M. \label{3.23}
       \eeq

\begin{example}
\label{e7}
 We consider two infinite numbers
 \[
 A=1\mbox{\ding{172}}^{18}(-5)\mbox{\ding{172}}^{2.4}
(-3)\mbox{\ding{172}}^{1}, \hspace{1cm} B=-1\mbox{\ding{172}}^{1}
0.7\mbox{\ding{172}}^{-3}
\]
 and calculate the
product  $C=B \cdot A$. The first partial product $C_1$ is   equal
to
\[
C_{1} = 0.7\mbox{\ding{172}}^{-3} \cdot A =
0.7\mbox{\ding{172}}^{-3}(\mbox{\ding{172}}^{18}-5\mbox{\ding{172}}^{2.4}-
3\mbox{\ding{172}}^{1})=
\]
\[
0.7\mbox{\ding{172}}^{15}-3.5\mbox{\ding{172}}^{-0.6}-2.1\mbox{\ding{172}}^{-2}
=
0.7\mbox{\ding{172}}^{15}(-3.5)\mbox{\ding{172}}^{-0.6}(-2.1)\mbox{\ding{172}}^{-2}.
\]
The second partial product, $C_2$, is computed analogously
\[
C_{2} = -\mbox{\ding{172}}^{1} \cdot A =
-\mbox{\ding{172}}^{1}(\mbox{\ding{172}}^{18}-5\mbox{\ding{172}}^{2.4}
-3\mbox{\ding{172}}^{1})=
-\mbox{\ding{172}}^{19}5\mbox{\ding{172}}^{3.4}3\mbox{\ding{172}}^{2}.
\]
Finally,   the product $C$ is equal to
\[
\begin{tabular}{cr}\hspace {11mm}$C =  C_{1} + C_{2}   =
 -1\mbox{\ding{172}}^{19}0.7\mbox{\ding{172}}^{15}
5\mbox{\ding{172}}^{3.4}3\mbox{\ding{172}}^{2}
(-3.5)\mbox{\ding{172}}^{-0.6}(-2.1)\mbox{\ding{172}}^{-2}.$ &
\end{tabular} \hspace{5mm}
 \Box
\]
%Execution of multiplication with the numbers $A$ and $B$ using the
%Infinity Calculator is illustrated in Fig.~\ref{Big_paper0}. Note
%that since the left and the right operands are written in their
%own windows, parenthesis for negative grossdigits are omitted.

 \end{example}

In the operation of \textit{division} of a   number $C$ by a
number $B$ from (\ref{3.20}), we obtain a result $A$ and a
reminder $R$ (that can be also equal to zero), i.e., $C =A \cdot
B+R$. The number $A$ is constructed as follows.  The first
grossdigit $a_{k_{K}}$ and the corresponding maximal exponent
${k_{K}}$ are established from the equalities
 \beq
a_{k_{K}}=c_{l_{L}}/ b_{m_{M}}, \hspace{4mm}  k_{K} = l_{L}-
m_{M}.
 \label{3.25}
       \eeq
Then the first partial reminder  $R_1$ is calculated as
 \beq
R_1= C - a_{k_{K}}\mbox{\ding{172}}^{k_{K}} \cdot B.
\label{3.25.0}
       \eeq
If $R_1 \neq 0$ then the number $C$ is substituted by $R_1$ and
the process is repeated with a complete analogy.  The grossdigit
$a_{k_{K-i}}$, the corresponding grosspower   $k_{K-i}$ and the
partial reminder $R_{i+1}$ are computed by  formulae
(\ref{3.25.1}) and (\ref{3.25.2}) obtained from (\ref{3.25}) and
(\ref{3.25.0}) as follows: $l_{L}$ and $c_{l_{L}}$ are substituted
by the highest grosspower $n_i$ and the corresponding grossdigit
$r_{n_{i}}$ of the partial reminder $R_{i}$ that, in  turn,
substitutes $C$:
        \beq
a_{k_{K-i}}=r_{n_i}/ b_{m_{M}}, \hspace{4mm}  k_{K-i} = n_i-
m_{M}.
 \label{3.25.1}
       \eeq
 \beq
R_{i+1}= R_{i} - a_{k_{K-i}}\mbox{\ding{172}}^{k_{K-i}}\cdot B,
\hspace{4mm}
  i \ge  1. \label{3.25.2}
       \eeq
The process stops when a partial reminder equal to zero is found
(this means that the final reminder $R=0$) or when a required
accuracy of the result is reached.

\begin{example}
\label{e8} Let us divide the number
$C=-10\mbox{\ding{172}}^{3}16\mbox{\ding{172}}^{0}42\mbox{\ding{172}}^{-3}$
by the number $B=5\mbox{\ding{172}}^{3}7$. For these numbers we
have
\[
l_{L} = 3, \hspace{2mm} m_{M}= 3, \hspace{2mm} c_{l_{L}}= -10,
\hspace{2mm} b_{m_{M}}= 5.
\]
It follows immediately from (\ref{3.25}) that
$a_{k_K}\mbox{\ding{172}}^{k_K}=-2\mbox{\ding{172}}^{0}$.  The
first partial reminder  $R_1$ is calculated as
\[
R_1=
-10\mbox{\ding{172}}^{3}16\mbox{\ding{172}}^{0}42\mbox{\ding{172}}^{-3}
- (-2\mbox{\ding{172}}^{0}) \cdot 5\mbox{\ding{172}}^{3}7=
\]
\[
-10\mbox{\ding{172}}^{3}16\mbox{\ding{172}}^{0}42\mbox{\ding{172}}^{-3}
+10\mbox{\ding{172}}^{3}14\mbox{\ding{172}}^{0} =
30\mbox{\ding{172}}^{0}42\mbox{\ding{172}}^{-3}.
\]
By a complete analogy we should construct
$a_{k_{K-1}}\mbox{\ding{172}}^{k_{K-1}}$   by rewriting
(\ref{3.25}) for $R_1$. By doing so we obtain equalities
\[
30=a_{k_{K-1}}\cdot 5, \hspace{4mm} 0 = k_{K-1} + 3
\]
and, as the result,  $a_{k_{K-1}}\mbox{\ding{172}}^{k_{K-1}}=
6\mbox{\ding{172}}^{-3}$. The second partial reminder   is
\[
R_2= R_1 - 6\mbox{\ding{172}}^{-3} \cdot 5\mbox{\ding{172}}^{3}7=
 30\mbox{\ding{172}}^{0}42\mbox{\ding{172}}^{-3} -
 30\mbox{\ding{172}}^{0}42\mbox{\ding{172}}^{-3} = 0.
\]
Thus, we can conclude that the reminder $R=R_2= 0$ and the final
result of division is
$A=-2\mbox{\ding{172}}^{0}6\mbox{\ding{172}}^{-3}$.

Let us now substitute the grossdigit  42 by 40 in $C$ and   divide
this new number
$\widetilde{C}=-10\mbox{\ding{172}}^{3}16\mbox{\ding{172}}^{0}40\mbox{\ding{172}}^{-3}$
by the same number $B=5\mbox{\ding{172}}^{3}7$. This operation
gives us the same result
$\widetilde{A}_2=A=-2\mbox{\ding{172}}^{0}6\mbox{\ding{172}}^{-3}$
(where subscript 2 indicates that two partial
reminders\index{partial reminder} have been obtained) but with the
reminder $\widetilde{R}=\widetilde{R}_2=
-2\mbox{\ding{172}}^{-3}$. Thus, we obtain $\widetilde{C} = B
\cdot \widetilde{A}_2+\widetilde{R}_2$. If we want to continue the
procedure of division, we obtain
$\widetilde{A}_3=-2\mbox{\ding{172}}^{0}6\mbox{\ding{172}}^{-3}(-0.4)\mbox{\ding{172}}^{-6}$
with the reminder $\widetilde{R}_3= 0.28\mbox{\ding{172}}^{-6}$.
Naturally, it follows $\widetilde{C} = B \cdot
\widetilde{A}_3+\widetilde{R}_3$. The process continues until a
partial reminder  $\widetilde{R}_i = 0$ is found   or when a
required accuracy of the result will be reached. \hfill $\Box$
 \end{example}

\section{Examples of   problems where   computations with new numerals can be useful}
\label{s4}

\subsection{The work with infinite sequences}
\label{s4.2}

We start by reminding   traditional definitions of the infinite
sequences and subsequences.  An \textit{infinite sequence}
$\{a_n\}, a_n \in A, n \in \mathbb{N},$ is a function having as
the domain the set of natural numbers, $\mathbb{N}$, and as the
codomain  a set $A$. A \textit{subsequence} is   a sequence from
which some of its elements have been removed. The IUA allows us to
prove the following result.
\begin{theorem}
\label{t2} The number of elements of any infinite sequence is less
or equal to~\ding{172}.
\end{theorem}

\textit{Proof.}  The IUA  states that the set $\mathbb{N}$ has
\ding{172} elements. Thus, due to the sequence definition given
above, any sequence having $\mathbb{N}$ as the domain  has
\ding{172} elements.

The notion of subsequence is introduced as a sequence from which
some of its elements have been removed. Thus, this definition
gives infinite sequences having the number of members less than
grossone.  \hfill $\Box$

One of the immediate consequences of the understanding of this
result is that any sequential process can have at maximum
\ding{172} elements. Due to Postulate 1, it depends on the chosen
numeral system which numbers among  \ding{172} members of the
process we can observe.

\begin{example}
\label{e12} Let us consider the set, $\widehat{\mathbb{N}}$, of
extended natural numbers from (\ref{4.2.2}). Then, starting from
the number 1, the process of the sequential counting can arrive at
maximum to \ding{172}
\[
\underbrace{1,2,3,4,\hspace{1mm}  \ldots \hspace{1mm}
\mbox{\ding{172}}-2,\hspace{1mm}
 \mbox{\ding{172}}-1,
\mbox{\ding{172}}}_{\mbox{\ding{172}}},  \mbox{\ding{172}}+1,
\mbox{\ding{172}}+2, \mbox{\ding{172}}+3, \ldots
\]
 Starting from 3 it   arrives at maximum
to $\mbox{\ding{172}}+2$
\[
\begin{tabular}{cr}\hspace {20mm}$1,2,\underbrace{3,4,\hspace{1mm}  \ldots \hspace{1mm}
\mbox{\ding{172}}-2,\hspace{1mm}
 \mbox{\ding{172}}-1,
\mbox{\ding{172}},  \mbox{\ding{172}}+1,
\mbox{\ding{172}}+2}_{\mbox{\ding{172}}}, \mbox{\ding{172}}+3,
  \ldots$ &
\hspace {13mm}
 $\Box$
\end{tabular}
\]
 \end{example}

It becomes appropriate now to define the \textit{complete
sequence} as an infinite sequence  containing \ding{172} elements.
For example, the sequence   of natural numbers  is complete, the
sequences of even  and odd natural numbers  are not complete.
Thus, the IUA imposes a more precise description of infinite
sequences. To define a sequence $\{a_n\}$ it is not sufficient
just to give a formula for~$a_n$, we should determine (as it
happens for sequences having a finite number of elements) the
first and the last elements of the sequence. If the number of the
first element is equal to one, we can use the record $\{a_n: k \}$
where $a_n$ is, as usual, the general element of the sequence and
$k$ is the number (that can be finite or infinite) of members of
the sequence.

\begin{example}
\label{e13} Let us consider the following two sequences, $\{a_n\}$
and $\{c_n\}$:
\[
 \{a_n\} =  \{ 5,\hspace{3mm} 10,\hspace{3mm} \ldots \hspace{3mm} 5(\mbox{\ding{172}}-1),\hspace{3mm}
5\mbox{\ding{172}} \},
\]
\beq \{b_n\}  = \{  5,\hspace{3mm}10,\hspace{3mm} \ldots
\hspace{3mm} 5 (\frac{2\mbox{\ding{172}}}{5}-1),\hspace{3mm}
5\cdot \frac{2\mbox{\ding{172}}}{5} \},
 \label{3.7.1}
       \eeq
\beq
 \{c_n\}  = \{ 5,\hspace{3mm} 10,\hspace{3mm} \ldots \hspace{3mm}
5 (\frac{4\mbox{\ding{172}}}{5}-1),\hspace{3mm} 5\cdot
\frac{4\mbox{\ding{172}}}{5} \}.
 \label{3.7.2}
       \eeq
They  have the same general element    $a_n=b_n=c_n=5n$ but they
are different because they have different numbers of members.  The
first sequence has \ding{172} elements and is thus complete,  the
other two sequences are not complete: $\{b_n\}$ has
$\frac{2\mbox{\ding{172}}}{5}$ elements  and $\{c_n\}$ has
$\frac{4\mbox{\ding{172}}}{5}$ members.   \hfill
 $\Box$
 \end{example}

In connection with this definition the following natural question
arises inevitably. Suppose that we have two sequences, for
example, $\{b_n:\frac{2\mbox{\ding{172}}}{5}\}$ and
$\{c_n:\frac{4\mbox{\ding{172}}}{5}\}$ from (\ref{3.7.1}) and
(\ref{3.7.2}). Can we create a new sequence, $\{d_n:k\}$, composed
from both of them, for instance, as it is shown below
 \[
b_1,\hspace{1mm} b_2,\hspace{1mm} \ldots \hspace{1mm}
b_{\frac{2\mbox{\tiny{\ding{172}}}}{5}-2},\hspace{1mm}
b_{\frac{2\mbox{\tiny{\ding{172}}}}{5}-1},\hspace{1mm}b_{\frac{2\mbox{\tiny{\ding{172}}}}{5}},\hspace{1mm}
c_1,\hspace{1mm} c_2,\hspace{1mm} \ldots \hspace{1mm}
c_{\frac{4\mbox{\tiny{\ding{172}}}}{5}-2},\hspace{1mm}
c_{\frac{4\mbox{\tiny{\ding{172}}}}{5}-1},\hspace{1mm}
c_{\frac{4\mbox{\tiny{\ding{172}}}}{5}}
 \]
and which will be the value of the number of its elements $k$?

The answer is `no' because  due to the  definition of the infinite
sequence, a sequence can be at maximum complete,  i.e., it cannot
have more than $\mbox{\ding{172}}$  elements. Starting from the
element $b_1$ we can arrive at maximum to the element
$c_{\frac{3\mbox{\tiny{\ding{172}}}}{5}}$  being the element
number \ding{172} in the sequence $\{d_n:k\}$ which we try to
construct. Therefore, $k=\mbox{\ding{172}}$ and
\[
\underbrace{b_1,\hspace{1mm}  \ldots \hspace{1mm}
b_{\frac{2\mbox{\tiny{\ding{172}}}}{5}},\hspace{1mm}
 c_1,\hspace{1mm}  \ldots
c_{\frac{3\mbox{\tiny{\ding{172}}}}{5}}}_{\mbox{\ding{172}
elements}}, \hspace{1mm}
\underbrace{c_{\frac{3\mbox{\tiny{\ding{172}}}}{5}+1}, \ldots
\hspace{1mm}
c_{\frac{4\mbox{\tiny{\ding{172}}}}{5}}}_{\frac{\mbox{\tiny{\ding{172}}}}{5}
 \mbox{ elements }}.
\]
The remaining members of the sequence
$\{c_n:\frac{4\mbox{\ding{172}}}{5}\}$ will form the second
sequence, $\{g_n: l \}$ having $l=
\frac{4\mbox{\ding{172}}}{5}-\frac{3\mbox{\ding{172}}}{5} =
\frac{\mbox{\ding{172}}}{5}$  elements. Thus, we have formed two
sequences, the first of them is complete and the second is not.

It is important to emphasize that the above consideration on the
infinite sequences allows us to deal with recursively defined
sets. Since such a set is constructed sequentially by a process,
it  can have at maximum \ding{172} elements.

To conclude this subsection, let us return to Hilbert's  paradox
of the Grand Hotel presented in Section~\ref{s1}. In the paradox,
the number of the rooms in the Hotel is countable. In our
terminology this means that it has \ding{172} rooms. When a new
guest arrives, it is proposed to move   the guest occupying room 1
to room 2, the guest occupying room 2 to room 3, etc. Under the
IUA this procedure does not help  because the guest from room
\ding{172} should be moved to room \ding{172}+1 and the Hotel has
only \ding{172} rooms. Thus, when the Hotel is full, no more new
guests can be accommodated -- the result corresponding perfectly
to Postulate~3 and the situation taking place in normal hotels
with a finite number of rooms.

\subsection{Calculating divergent series}
\label{s4.3}

Let us show how the new approach can be applied in  such an
important area as theory of divergent series. We consider two
infinite series $S_1=10+10+10+\ldots$ and $S_2=3+3+3+\ldots$  The
traditional analysis gives us a very poor answer that both of them
diverge to infinity. Such operations as, e.g., $\frac{S_2}{S_1}$
and $S_2 - S_1$
 are not defined.

Now, when we are able to express not only different finite numbers
but also different infinite numbers, it is necessary to indicate
explicitly the number of items in the sums $S_1$ and $S_2$ and it
is not important if it is finite or infinite. To calculate the sum
it is necessary that  the number of items and the result are
expressible in the numeral system used for calculations. It is
important to notice that even though a sequence cannot have more
than \ding{172} elements, the number of items in a series can be
greater than grossone because the process of summing up is not
necessary  executed by a sequential adding items.

Let us suppose that the   series $S_1$ has $k$ items and $S_2$ has
$n$ items. We can then define sums (that can have a finite or an
infinite number of items),
$$S_1(k)=\underbrace{10+10+10+\ldots+10}_k, \hspace{1cm} S_2(n)=\underbrace{3+3+3+\ldots+3}_n,$$
 calculate them, and execute
arithmetical operations with the obtained results. The sums then
are obviously calculated as $S_1(k)=10k$ and $S_2(n)=3n$.  If, for
instance, $k=n=5\mbox{\ding{172}}$ then we obtain
$S_1(5\mbox{\ding{172}})=50\mbox{\ding{172}}$,
$S_2(5\mbox{\ding{172}})=15\mbox{\ding{172}}$ and
\[
S_2(5\mbox{\ding{172}}) /  S_1(5\mbox{\ding{172}}) = 0.3.
\]
Analogously, if  $k=3\mbox{\ding{172}}$ and
$n=10\mbox{\ding{172}}$ we obtain
$S_1(3\mbox{\ding{172}})=30\mbox{\ding{172}}$,
$S_2(\mbox{\ding{172}})=30\mbox{\ding{172}}$ and  it follows
 $S_2(\mbox{\ding{172}}) -
S_1(3\mbox{\ding{172}})=0$.

If   $k=3\mbox{\ding{172}}4$ (we remind that we use here a shorter
way to write down this infinite number, the complete record is
$3\mbox{\ding{172}}^{1}4\mbox{\ding{172}}^{0}$) and
$n=10\mbox{\ding{172}}$ we obtain
$S_1(3\mbox{\ding{172}}4)=30\mbox{\ding{172}}40$,
$S_2(\mbox{\ding{172}})=30\mbox{\ding{172}}$ and  it follows
 \[
S_1(3\mbox{\ding{172}}4) -  S_2(\mbox{\ding{172}}) =
30\mbox{\ding{172}}40 - 30\mbox{\ding{172}}=   40.
 \]
 \[
S_1(3\mbox{\ding{172}}2) / S_2(\mbox{\ding{172}}) =
30\mbox{\ding{172}}20 / 30\mbox{\ding{172}}=
1\mbox{\ding{172}}^{0}0.66667\mbox{\ding{172}}^{-1} > 0.
 \]

We conclude this subsection by studying the series
$\sum_{i=1}^{\infty}\frac{1}{2^i}$. It is known that it converges
to one. However, we are able to give a more precise answer. Due to
Postulate~3, the formula
$$\sum_{i=1}^{k}\frac{1}{2^i}=1-\frac{1}{2^k}$$
can be used directly for   infinite $k$, too. For example, if
$k=\mbox{\ding{172}}$ then
$$\sum_{i=1}^{\mbox{\small{\ding{172}}}}\frac{1}{2^i}=1-\frac{1}{2^{\mbox{\tiny{\ding{172}}}}}$$
where $\frac{1}{2^{\mbox{\tiny{\ding{172}}}}}$ is infinitesimal.
Thus, the traditional answer $\sum_{i=1}^{\infty}\frac{1}{2^i}=1$
is a finite approximation to our more precise result using
infinitesimals. More examples related to series can be found in
\cite{chaos}.

\subsection{Calculating limits and expressing irrational numbers}
\label{s4.5}

 Let us now discuss the problem of calculation of limits from the point of
view of our approach. In traditional analysis, if a limit $\lim_{x
\rightarrow a}f(x)$ exists, then it gives us a very poor -- just
one  value -- information about the behavior of $f(x)$ when $x$
tends to $a$. Now we can obtain significantly   richer information
because we are able to calculate $f(x)$ directly at any finite,
infinite, or infinitesimal point that can be expressed by the new
positional system even if the limit does not exist.

Thus, limits  equal to infinity can be substituted by  precise
infinite numerals  and limits equal to zero can be substituted by
precise infinitesimal numerals\footnote{Naturally, if we speak
about limits of sequences, $\lim_{n \rightarrow \infty}a(n)$, then
$ n \in \mathbb{N}$ and, as a consequence, it follows that $n$
should be less than or equal to grossone.}.  This is very
important for practical computations because these substitutions
eliminate  indeterminate forms.

\begin{example}
\label{e15} Let us consider the following two limits
 \[
\lim_{x \rightarrow +\infty}(5x^3-x^2+10^{61})=
 +\infty, \hspace{1cm} \lim_{x \rightarrow +\infty}(5x^3-x^2)=
 +\infty.
 \]
Both give  us   the same result, $+\infty$, and it is not possible
to execute the operation
 \[
 \lim_{x \rightarrow +\infty}(5x^3-x^2+10^{61}) - \lim_{x \rightarrow +\infty}(5x^3-x^2).
 \]
that is an   indeterminate form of the type $\infty-\infty$ in
spite of the fact that for any finite $x$ it follows
 \beq
5x^3-x^2+10^{61} - (5x^3-x^2) = 10^{61}. \label{4.4.3}
 \eeq
The new approach allows us to calculate exact values of both
expressions, $5x^3-x^2+10^{61}$ and $5x^3-x^2$, at any infinite
(and infinitesimal) $x$ expressible in the chosen numeral system.
For instance, the choice $x=3\mbox{\ding{172}}^{2}$ gives the
value
\[
5(3\mbox{\ding{172}}^{2})^{3}-(3\mbox{\ding{172}}^{2})^{2}+10^{61}=
135\mbox{\ding{172}}^{6}\mbox{\small-}9\mbox{\ding{172}}^{4}10^{61}
\]
for the first expression and
$135\mbox{\ding{172}}^{6}\mbox{\small-}9\mbox{\ding{172}}^{4}$ for
the second one. We can easily calculate the difference of these
two infinite numbers, thus obtaining the same result as we had for
finite values of $x$ in  (\ref{4.4.3}):
 \[
\begin{tabular}{cr}\hspace {32mm}$135\mbox{\ding{172}}^{6}\mbox{\small-}9\mbox{\ding{172}}^{4}10^{61} -
(135\mbox{\ding{172}}^{6}\mbox{\small-}9\mbox{\ding{172}}^{4}) =
10^{61}.$ & \hspace {21mm}
 $\Box$
\end{tabular}
\]
\end{example}

It is necessary to emphasize the fact  that expressions   can be
calculated even when their limits do not exist. Thus, we obtain a
very powerful tool for studying divergent processes.
\begin{example}
\label{e16}
 The   limit
$ \lim_{n \rightarrow +\infty}f(n),$ $f(n)=(-1)^n n^3$, does not
exist. However, we can easily calculate expression $(-1)^n n^3$ at
different infinite points $n$. For instance, for
$n=\mbox{\ding{172}}$ it follows
$f(\mbox{\ding{172}})=\mbox{\ding{172}}^3$ because grossone is
even and for the odd $n=0.5\mbox{\ding{172}}-1$ it follows
\[
\hspace{18mm}
f(0.5\mbox{\ding{172}}-1)=-(0.5\mbox{\ding{172}}-1)^3=
-0.125\mbox{\ding{172}}^{3}0.75\mbox{\ding{172}}^{2}\mbox{\small-}1.5\mbox{\ding{172}}^{1}1.
\hspace{15mm} \Box
\]
\end{example}

Limits with  the argument tending to zero can be considered
analogously. In this case, we can calculate the corresponding
expression at any infinitesimal point using the new positional
system and   obtain a significantly more rich information.

\begin{example}
\label{e17} If $x$ is a fixed finite number then
 \beq
 \lim_{h \rightarrow 0}\frac{(x+h)^2-x^2}{h}= 2x.
\label{4.5}
 \eeq
In the new positional system  we obtain
 \beq
\frac{(x+h)^2-x^2}{h}= 2x + h.
 \label{4.6}
 \eeq
If, for instance,   $h=\mbox{\ding{172}}^{-1}$, the answer is
$2x\mbox{\ding{172}}^{0}\mbox{\ding{172}}^{-1}$, if
$h=4.2\mbox{\ding{172}}^{-2}$ we obtain the value
$2x\mbox{\ding{172}}^{0}4.2\mbox{\ding{172}}^{-2}$, etc. Thus, the
value of  the limit (\ref{4.5}), for a finite $x$, is just the
finite approximation of the number (\ref{4.6}) having finite and
infinitesimal parts. \hfill $\Box$
\end{example}

Let us make a remark regarding irrational numbers. Among their
properties, they are characterized by the fact that we do not know
any numeral system that would allow  us to express them by a
finite number of symbols used to express other numbers. Thus,
special numerals ($e, \pi, \sqrt{2}, \sqrt{3} $, etc.) are
introduced by describing their properties in a way (similarly, all
other numerals, e.g., symbols `0' or `1', are introduced also by
describing their properties). These special symbols are then used
in analytical transformations together with ordinary numerals.

For example,   it is possible to work directly with the symbol $e$
in analytical transformations by applying suitable rules defining
this number together with numerals taking part in a chosen numeral
system $\mathcal{S}$. At the end of transformations, the obtained
result will be be expressed in numerals from $\mathcal{S}$ and,
probably, in terms of $e$. If it is then required to execute some
\textit{numerical} computations, this means that it is necessary
to substitute $e$ by a numeral (or numerals) from $\mathcal{S}$
that will allow us to approximate  $e$ in some way.

The same situation takes place when one uses the new numeral
system, i.e., while we work analytically we use just the symbol
$e$ in our expressions and then, if we wish to work numerically we
should pass to approximations. The new numeral system opens a new
perspective on the problem of the expression of irrational
numbers. Let us consider one of the possible ways to obtain an
approximation of $e$, i.e., by using the limit
 \beq
e = \lim_{n \rightarrow +\infty}(1+\frac{1}{n})^n  =
2.71828182845904\ldots
 \label{3.200}
       \eeq
In our numeral system the expression $(1+\frac{1}{n})^n$ can be
written directly for finite and/or infinite values of $n$. For
$n=\mbox{\ding{172}}$ we obtain the number $e_0$ designated so in
order to distinguish it from the record (\ref{3.200})
 \beq
e_0 = (1+\frac{1}{\mbox{\ding{172}}})^{\mbox{\tiny{\ding{172}}}}=
(\mbox{\ding{172}}^{0}\mbox{\ding{172}}^{-1})^{\mbox{\tiny{\ding{172}}}}.
\label{3.201}
       \eeq
It becomes clear from this record why the number $e$ cannot be
expressed in a positional numeral system with a finite base. Due
to the definition of a sequence under the IUA,  such a system can
have at maximum \ding{172} numerals -- digits -- to express
fractional part of a number (see section~\ref{s4.4} for details)
and, as it can be seen from (\ref{3.201}), this quantity is not
sufficient for $e$ because the item
$\frac{1}{\mbox{\ding{172}}^{\mbox{\tiny{\ding{172}}}}}$ is
present in it.

Naturally, it is also possible to construct more exotic   $e$-type
numbers by substituting \ding{172} in (\ref{3.201}) by any
infinite number written in the new positional system  with
infinite base. For example, if we substitute  \ding{172} in
(\ref{3.201}) by $\mbox{\ding{172}}^{2}$ we obtain the number
\[
e_1 =
(1+\frac{1}{\mbox{\ding{172}}^{2}})^{\mbox{\tiny{\ding{172}}}^{2}}=
(\mbox{\ding{172}}^{0}\mbox{\ding{172}}^{-2})^{\mbox{\tiny{\ding{172}}}^{2}}.
\]
The numbers considered above take their origins in the limit
(\ref{3.200}). Similarly, other formulae leading to approximations
of $e$ expressed in traditional numeral systems give us other new
numbers that can be expressed in the new numeral system. The same
way of reasoning can be used with respect to other irrational
numbers, too.

\subsection{Measuring infinite sets with elements defined by formulae}
\label{s4.1}

We have already discussed in Section~\ref{s2} how   we calculate
the number of elements for sets being results of the usual
operations (intersection, union, etc.) with finite sets and
infinite sets of the type $\mathbb{N}_{k,n}$. In order to have a
possibility to work with infinite sets having a more general
structure than the sets $\mathbb{N}_{k,n}$, we need to develop
more powerful instruments. Suppose that we have an integer
function $g(i)
> 0$ strictly increasing on indexes $i=1,2,3, \ldots$ and we wish to
know how many elements are there in the set
$$G = \{ g(1), g(2), g(3), \ldots \}.$$  In our terminology this
question has no any sense because of the following reason.

In the finite case, to define  a   set it is not sufficient to say
that it is  finite.   It is necessary to indicate its number  of
elements explicitly as, e.g., in this example
\[
G_1 = \{ g(i): 1 \le i \le 5 \},
 \]
or implicitly, as it is made here:
 \beq
  G_2 = \{ g(i): i
\ge 1, \,\,\,\, 0 < f(i)\le b \}, \label{4.3}
 \eeq
 where   $b$ is   finite.

Now we have mathematical tools to indicate the number of elements
for   infinite sets, too. Thus, analogously to the finite case and
due to Postulate~3, it is not sufficient to say that a set has
infinitely many elements. It is necessary to indicate its number
of elements explicitly or implicitly. For instance, the number of
elements of the set
\[
G_3 = \{ g(i): 1 \le i \le   \mbox{\ding{172}}^{10}   \}
\]
is indicated explicitly: the set $G_3$ has
$\mbox{\ding{172}}^{10}$ elements.

If a set is given in the form (\ref{4.3}) where $b$ is infinite,
then its number of elements, $J$, can be determined  as
 \beq
 J = \max
\{i : g(i) \le b \} \label{4.3.0}
 \eeq
 if we are able to determine the inverse
function $g^{-1}(x)$ for $g(x)$. Then, $J= [ g^{-1}(b)]$, where
$[u]$ is integer part of $u$. Note that if $b=\mbox{\ding{172}}$,
then the set $G_2 \subseteq \mathbb{N}$ since all its elements are
integer, positive, and $g(i) \le \mbox{\ding{172}}$ due to
(\ref{4.3.0}).

\begin{example}
\label{e9} Let us consider the following set, $A_1(k,n)$, having
$g(i) = k+n(i-1),$
\[
A_1(k,n) = \{ g(i) : i \ge 1, \,\,\,\,g(i) \le \mbox{\ding{172}}
\}, \hspace{1cm} 1 \le k \le n, \hspace{3mm} n \in \mathbb{N}.
\]
It follows   from the IUA that $A_1(k,n) = \mathbb{N}_{k,n}$ from
(\ref{3.3}).   By applying (\ref{4.3.0}) we find for $A_1(k,n)$
its number of elements
\[
\begin{tabular}{cr}\hspace {11mm}$J_1(k,n)= [ \frac{\mbox{\ding{172}}-k}{n}+1]= [
\frac{\mbox{\ding{172}}-k}{n}]+1= \frac{\mbox{\ding{172}}}{n}-1+1=
\frac{\mbox{\ding{172}}}{n}.$ & \hspace {23mm}
 $\Box$
\end{tabular}
\]
 \end{example}

\begin{example}
\label{e10} Analogously, the  set
\[
A_2(k,n,j) = \{ k+ni^j : i \ge 0,\,\,\,\, 0 < k+ni^j \le
\mbox{\ding{172}} \}, \hspace{7mm} 0 \le k < n, \hspace{2mm}  n
\in \mathbb{N}, \hspace{2mm} j \in \mathbb{N},
\]
has $J_2(k,n,j)=  [ \sqrt[j]{\frac{\mbox{\ding{172}}-k}{n}}]$
elements. \hfill $\Box$
 \end{example}

\subsection{Measuring infinite sets of numerals and their comparison}
\label{s4.4}

Let us calculate the number of elements in some well known
infinite sets of numerals using the designation $|A|$ to indicate
the number of elements of a set $A$.

\begin{theorem}
\label{t3}The number of elements  of the  set, $\mathbb{Z}$, of
integers is $|\mathbb{Z}|= 2\mbox{\ding{172}}1$.
\end{theorem}

\textit{Proof.} The set $\mathbb{Z}$  contains \ding{172} positive
numbers, \ding{172} negative numbers, and zero. Thus,
 \beq
|\mathbb{Z}|= \mbox{\ding{172}} + \mbox{\ding{172}} + 1 =
2\mbox{\ding{172}}1. \hspace{5mm}
 \Box
 \label{3.102}
       \eeq

Traditionally,   rational numbers are defined as ratio of two
integer numbers. The new approach allows us to calculate the
number of numerals in a fixed numeral system. Let us consider a
numeral system $\mathbb{Q}_1$ containing numerals of the form
 \beq
 \frac{p}{q}, \hspace{5mm} p   \in \mathbb{Z},\hspace{3mm} q \in \mathbb{Z},\,\,\, q \neq 0.
 \label{3.102.0}
       \eeq
\begin{theorem}
\label{t4.0}
 The number of elements   of the set, $\mathbb{Q}_1$,
of rational numerals of the type (\ref{3.102.0}) is
$|\mathbb{Q}_1|=4\mbox{\ding{172}}^{2}2\mbox{\ding{172}}^1$.
\end{theorem}

\textit{Proof.} It follows from Theorem~\ref{t3} that the
numerator of (\ref{3.102.0}) can be filled in by
$2\mbox{\ding{172}}1$ and the denominator by $2\mbox{\ding{172}}$
numbers.  Thus, number of all possible combinations is
\[
|\mathbb{Q}_1|= 2\mbox{\ding{172}}1 \cdot 2\mbox{\ding{172}} =
4\mbox{\ding{172}}^{2}2\mbox{\ding{172}}^1.
 \]

 \hfill  $\Box$

It is necessary to notice that in   Theorem~\ref{t4.0} we have
calculated different numerals and not different numbers. For
example, in the numeral system $\mathbb{Q}_1$ the number 0 can be
expressed by $2\mbox{\ding{172}}$ different numerals
\[
\frac{0}{-\mbox{\ding{172}}}, \,\, \frac{0}{-\mbox{\ding{172}+1}},
\,\, \frac{0}{-\mbox{\ding{172}+2}}, \ldots \frac{0}{-2}, \,\,
\frac{0}{-1}, \,\, \frac{0}{1},  \,\, \frac{0}{2}, \ldots
\frac{0}{\mbox{\ding{172}-2}}, \,\, \frac{0}{\mbox{\ding{172}-1}},
\,\, \frac{0}{\mbox{\ding{172}}}
\]
and numerals such as $\frac{-1}{-2}$ and $\frac{1}{2}$ have been
calculated as two different numerals. The following theorem
determines the number of elements of the set $\mathbb{Q}_2$
containing numerals of the form
 \beq
  -\frac{p}{q}, \,\, \frac{p}{q}, \hspace{5mm} p   \in \mathbb{N},\hspace{3mm} q \in
 \mathbb{N},
 \label{3.102.1}
       \eeq
and zero is represented by one symbol 0.

\begin{theorem}
\label{t4}The number of elements   of the set, $\mathbb{Q}_2$, of
rational numerals   of the type (\ref{3.102.1}) is
$|\mathbb{Q}_2|=2\mbox{\ding{172}}^{2}1$.
\end{theorem}

\textit{Proof.} Let us consider positive rational numerals.  The
form of the rational numeral $\frac{p}{q}$, the fact that $p, \, q
\in \mathbb{N},$ and the IUA impose that both $p$ and $q$ can
assume values from 1 to \ding{172}. Thus, the number of all
possible combinations is $\mbox{\ding{172}}^{2}$. The same number
of combinations we obtain for negative rational numbers and one is
added because we count zero as well. \hfill
 $\Box$

Let us now calculate the number of elements of the set,
$\mathbb{R}_{b}$, of real numbers expressed by numerals in the
positional system by the record\vspace*{-2mm}
 \beq
 (a_{n-1}a_{n-2} \ldots a_1 a_0 .
 a_{-1} a_{-2}  \ldots   a_{-(q-1)}     a_{-q})_b
   \label{3.103}
       \eeq
where the symbol $b$ indicates the radix\index{radix} of the
record and $n, \, q \in \mathbb{N}$.

\begin{theorem}
\label{t5} The number of elements   of the set, $\mathbb{R}_{b}$,
of numerals (\ref{3.103}) is $|\mathbb{R}_{b}| =
b^{2\mbox{\ding{172}}}$.
\end{theorem}

\textit{Proof.} In formula (\ref{3.103}) defining the type of
numerals we deal with there are two sequences of digits: the first
one, $a_{n-1}a_{n-2} \ldots a_1 a_0$, is used to express the
integer part of the number and the second, $a_{-1} a_{-2} \ldots
a_{-(q-1)} a_{-q}$, for its fractional part. Due to definition of
sequence and the IUA, each of them can have at maximum \ding{172}
elements. Thus, it can be at maximum \ding{172} positions on the
left of the dot
  and, analogously,
\ding{172} positions on the right of the dot. Every position can
be filled in by one of the $b$ digits from the
alphabet\index{alphabet} $\{ 0, 1, \ldots , b-1 \}$. Thus, we have
$b^{\mbox{\ding{172}}}$ combinations to express the integer part
of the number and the same quantity to express its fractional
part. As a result,   the positional numeral system using the
numerals of the form (\ref{3.103})   can express
$b^{2\mbox{\ding{172}}}$ numbers. \hfill
 $\Box$

Note that the result of theorem~\ref{t5} does not consider the
practical situation  of writing down concrete numerals. Obviously,
the number of numerals of the type (\ref{3.103}) that can be
written in practice is finite and depends on the chosen numeral
system for writing   digits.

It is worthwhile to notice also that the traditional point of view
on real numbers tells that there exist real numbers that can be
represented in positional systems by two different infinite
sequences of digits. In contrast, under the IUA all the numerals
represent different numbers. In addition, minimal and maximal
numbers expressible in $\mathbb{R}_{b}$ can be explicitly
indicated.
\begin{example}
\label{e14} For instance, in the decimal positional
system\index{positional system} $\mathbb{R}_{10}$ the numerals
\[
1.\underbrace{999\ldots99}_{\mbox{\ding{172} digits}},
\,\,\,\,\,\,\,\, 2.\underbrace{000\ldots00}_{\mbox{\ding{172}
digits}}
\]
represent different numbers and their difference is equal to
\[
2.\underbrace{000\ldots00}_{\mbox{\ding{172}
digits}}-1.\underbrace{999\ldots9}_{\mbox{\ding{172} digits}}
=0.\underbrace{000\ldots01}_{\mbox{\ding{172} digits}}.
\]
Analogously the smallest and the largest numbers expressible in
$\mathbb{R}_{10}$ can be easily indicated. They are, respectively,
\[
\hspace{25mm}-\underbrace{999\ldots9}_{\mbox{\ding{172}
digits}}.\underbrace{999\ldots9}_{\mbox{\ding{172} digits}},
\hspace{10mm} \underbrace{999\ldots9}_{\mbox{\ding{172}
digits}}.\underbrace{999\ldots9}_{\mbox{\ding{172} digits}}.
\hspace{25mm}
 \Box
 \]

 \end{example}

\begin{theorem}
\label{t6} The sets $\mathbb{Z}, \mathbb{Q}_1,\mathbb{Q}_2,$ and $
\mathbb{R}_{b}$ are not   monoids under addition.
 \end{theorem}

\textit{Proof.} The proof is obvious and is so omitted. \hfill
 $\Box$

\section{Relations to   results of Georg Cantor}
\label{s6}

It is obligatory to say in this occasion that the results
presented above should be considered as a more precise analysis of
the situation discovered by the genius of Cantor. He has proved,
by using his famous diagonal argument, that the number of elements
of the set  $\mathbb{N}$ is less than the number of real numbers
at the interval $[0,1)$ \textit{without calculating the latter}.
To do this he expressed  real numbers in a positional numeral
system. We have shown that this number will be different depending
on the radix $b$ used in the positional system to express real
numbers. However, all of the obtained numbers,
$b^{2\mbox{\ding{172}}}$, are more than the number of elements of
the set of natural numbers, \ding{172}, and, therefore, the
diagonal argument maintains its force.

We can now calculate the number of points of the interval $[0,1)$,
of a line, and of the $N$-dimensional space. To do this we need a
definition of the term  \textit{point}\index{point} and
mathematical tools to indicate a point. Since this concept is one
of the most fundamental, it is very difficult to find an adequate
definition. If we accept (as   is usually done in   modern
Mathematics) that the \textit{point}  $x$ in an $N$-dimensional
space is determined by $N$ numerals called \textit{coordinates of
the point}
\[
(x_1,x_2, \ldots x_{N-1},x_{N}) \in \mathbb{S}^N,
\]
where $\mathbb{S}^N$ is a set of numerals,    then we can indicate
the point $x$ by its coordinates   and we are able to execute the
required calculations. It is worthwhile to emphasize  that we have
not postulated that $(x_1,x_2, \ldots $ $ x_{N-1},x_{N})$ belongs
to the $N$-dimensional set, $\mathbb{R}^N$, of real numbers  as it
is usually done because we can express coordinates only by
numerals and, as we have shown above, different choices of numeral
systems lead to various sets of numerals.

We should decide now which numerals we shall use to express
coordinates of the points. Different variants can be chosen
depending on the precision level we want to obtain. For example,
if the numbers $0 \le x < 1$ are expressed in the form
$\frac{p-1}{\mbox{\ding{172}}}, p  \in \mathbb{N}$, then the
smallest positive number we can distinguish is
$\frac{1}{\mbox{\ding{172}}}$. Therefore, the interval $[0,1)$
contains the  following \ding{172} points
\[
0, \,\,\,\, \frac{1}{\mbox{\ding{172}}}, \,\,\,\,
\frac{2}{\mbox{\ding{172}}}, \,\,\,\,  \ldots \,\,\,\,
\frac{\mbox{\ding{172}}-2}{\mbox{\ding{172}}}, \,\,\,\,
\frac{\mbox{\ding{172}}-1}{\mbox{\ding{172}}}.
\]
Then, due to the IUA  and the definition of sequence, there are
\ding{172} intervals of the form $[a-1,a), a \in \mathbb{N},$ on
the ray  $x \ge 0$. Hence, this ray contains
$\mbox{\ding{172}}^{2}$ points and the whole line consists of
$2\mbox{\ding{172}}^{2}$ points.

If we need a higher precision, within each interval
\[
[a-1+\frac{i-1}{\mbox{\ding{172}}},
a-1+\frac{i}{\mbox{\ding{172}}}), \hspace{3mm} a,i \in \mathbb{N},
\]
we can distinguish again \ding{172} points and the number of
points within each interval $[a-1,a), a \in \mathbb{N},$ will
become  equal to $\mbox{\ding{172}}^{2}$. Consequently, the number
of the points  on the line will be equal to
$2\mbox{\ding{172}}^{3}$.

This situation   is a direct consequence of Postulate~2 and is
typical for natural sciences where it is well known that
instruments influence the results of observations. It is similar
as to   work with a microscope: we decide the level of the
precision we need and obtain a result which is dependent on the
chosen level of accuracy. If we need a more precise or a more
rough answer, we change the lens of our microscope.

Continuing the analogy with the microscope, we can also decide to
change our microscope with a new one. In our terms this means to
change the numeral system with another one. For instance, instead
of the numerals considered above, we choose a positional numeral
system to calculate the number of points within the interval
$[0,1)$; then, as we have already seen before, we are able to
distinguish $b^{\mbox{\ding{172}}}$ points of the form
\[
(.a_{-1} a_{-2}  \ldots a_{-(\mbox{\ding{172}}-1)}
a_{-\mbox{\ding{172}}})_b
\]
on it. Since the line contains $2\mbox{\ding{172}}$ unit
intervals, the whole number of points of this type on the line is
equal to $2\mbox{\ding{172}}b^{\mbox{\ding{172}}}$.

In this example of counting, we have changed the tool to calculate
the number of points within each interval, but used the old way to
calculate the number of intervals, i.e., by natural numbers.  If
we are not interested in subdividing the line at intervals and
want to obtain the number of the points on the line directly by
using positional numerals of the type (\ref{3.103}) with possible
infinite $n$ and $q$, then we are able to distinguish at maximum
$b^{2\mbox{\ding{172}}}$ points on the line.

 \begin{figure}[t]
  \begin{center}
    \epsfig{ figure = 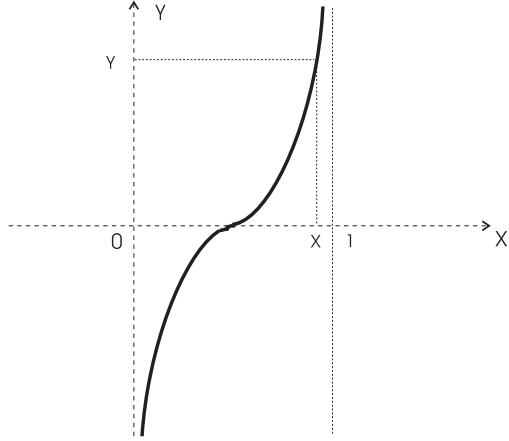, width = 2.6in, height = 2.3in,  silent = yes }
    \caption{Due to Cantor, the interval $(0,1)$ and the entire real number line have
the same number of points}
 \label{figura3}
  \end{center}
\end{figure}

 \begin{figure}[ht]
  \begin{center}
    \epsfig{ figure = 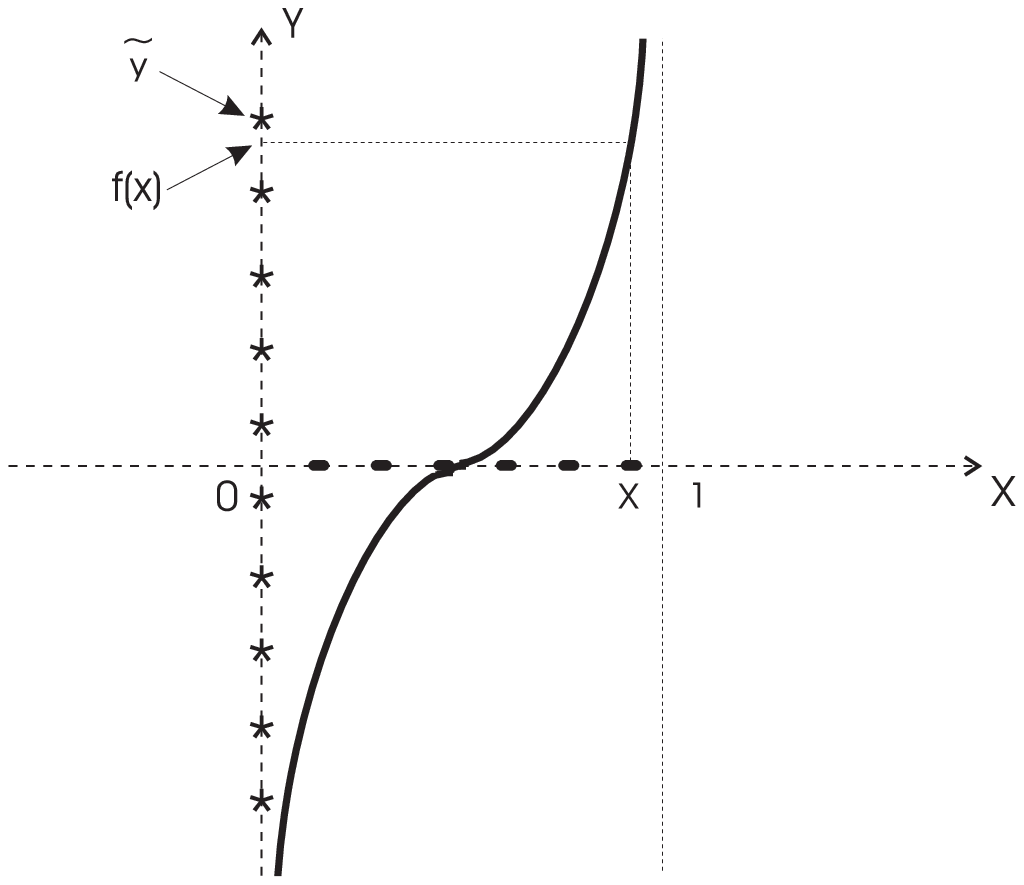, width = 2.6in, height = 2.3in,  silent = yes }
    \caption{Three independent mathematical objects:   the  set  $X_{\mathcal{S}_1}$ represented by dots,   the
    set
$Y_{\mathcal{S}_2}$ represented by stars, and   function
(\ref{4.4})}
 \label{Big_paper1}
  \end{center}
\end{figure}

Let us now return to the problem of comparison of infinite sets
and consider Cantor's famous result showing   that the number of
points over the interval $(0,1)$ is equal to the number of points
over the whole real line, i.e.,
 \beq
 |\mathbb{R}| = |(0,1)|.
\label{1.11}
 \eeq
 The proof of this counterintuitive fact is
given by establishing a one-to-one correspondence between the
elements of the two sets. Such a mapping can be done by using  for
example the function
 \beq
 y=\tan(0.5 \pi (2x-1)),\hspace{15mm}
x \in (0,1),
\label{4.4}
 \eeq
illustrated in Fig.~\ref{figura3}. Cantor shows by using
Fig.~\ref{figura3} that to any point $x \in (0,1)$ a point $y \in
(-\infty,\infty)$ can be associated and vice versa. Thus, he
concludes that
 the requested one-to-one correspondence  between the
sets $\mathbb{R}$ and $ (0,1)$ has been established and,
therefore, this proves (\ref{1.11}).

Our point of view is different: the number of elements is an
intrinsic characteristic of each set (for both finite and infinite
cases) that does not depend on any object outside the set. Thus,
in  Cantor's example from Fig.~\ref{figura3}  we have (see
Fig.~\ref{Big_paper1}) three mathematical objects: (i) a set,
$X_{\mathcal{S}_1}$, of points over the interval $ (0,1)$ which we
are able to distinguish using a numeral system $\mathcal{S}_1$;
(ii) a set, $Y_{\mathcal{S}_2}$, of points over the vertical real
line which we are able to distinguish using a numeral system
$\mathcal{S}_2$; (iii) the function (\ref{4.4}) described using a
numeral system $\mathcal{S}_3$. All these three mathematical
objects are independent  each other. The sets $X_{\mathcal{S}_1}$
and $Y_{\mathcal{S}_2}$  can have the same or different number of
elements.

Thus, we are not able to evaluate $f(x)$ at \textit{any} point
$x$. We are able to do this only at points from
$X_{\mathcal{S}_1}$.   Of course, in order to be able to execute
these evaluations it is necessary to conciliate the numeral
systems $\mathcal{S}_1,  \mathcal{S}_2,$ and $\mathcal{S}_3$. The
fact that we have made  evaluations of $f(x)$ and have obtained
the corresponding values does not influence minimally   the
numbers of elements of the sets $X_{\mathcal{S}_1}$ and
$Y_{\mathcal{S}_2}$. Moreover, it can happen  that the number
$y=f(x)$ cannot be expressed in the numeral system $\mathcal{S}_2$
and it is necessary to approximate it by a number $\widetilde{y}
\in \mathcal{S}_2$. This situation, very well known to computer
scientists, is represented in Fig.~\ref{Big_paper1}.

Let us remind one more famous example related to the one-to-one
correspondence and  taking its origins in studies of Galileo
Galilei: even numbers can be put in a one-to-one correspondence
with all natural numbers in spite of the fact that they are a part
of them:
 \beq
\begin{array}{lccccccc}
  \mbox{even numbers:}   & \hspace{5mm} 2, & 4, & 6, & 8,  & 10, & 12, & \ldots    \\

& \hspace{5mm} \updownarrow &  \updownarrow & \updownarrow  & \updownarrow  & \updownarrow  &  \updownarrow &   \\

  \mbox{natural numbers:}& \hspace{5mm}1, &  2, & 3, & 4 & 5,
       & 6,  &    \ldots \\
     \end{array}
\label{4.4.1}
 \eeq

Again, our view on this situation is different since we cannot
establish a one-to-one correspondence between the sets because
they are infinite and we, due to Postulate~1, are able to execute
only a finite number of operations. We cannot use the one-to-one
correspondence as an executable operation when it is necessary to
work with infinite sets.

However,   we   already know that the number of elements of the
set of natural numbers is equal to \ding{172} and \ding{172} is
even. Since the number of elements of the set of even numbers is
equal to $\frac{\mbox{\ding{172}}}{2}$, we   can write down not
only initial (as it is usually done traditionally) but also the
final part of (\ref{4.4.1})
  \beq
\begin{array}{cccccccccc}
 2, & 4, & 6, & 8,  & 10, & 12, & \ldots  &
\mbox{\ding{172}} -4,  &    \mbox{\ding{172}}  -2,   &    \mbox{\ding{172}}    \\
 \updownarrow &  \updownarrow & \updownarrow  &
\updownarrow  & \updownarrow  &  \updownarrow  & &
  \updownarrow    & \updownarrow   &
  \updownarrow
   \\
 1, &  2, & 3, & 4 & 5, & 6,   &   \ldots  &    \frac{\mbox{\ding{172}}}{2} - 2,  &
     \frac{\mbox{\ding{172}}}{2} - 1,  &    \frac{\mbox{\ding{172}}}{2}   \\
     \end{array}
\label{4.4.2}
 \eeq
concluding so (\ref{4.4.1})   in a complete accordance with
Postulate~3. Note that record (\ref{4.4.2}) does not affirms that
we have established the one-to-one correspondence among
\textit{all}  even numbers and a half of natural ones. We cannot
do this due to Postulate~1. The symbols `$\ldots$' indicate an
infinite number of numbers and we can execute only a finite number
of operations. However,  record (\ref{4.4.2}) affirms that for any
even number expressible in the chosen numeral system it is
possible to indicate the corresponding natural number in the lower
row of (\ref{4.4.2}).

We conclude the   paper  by the following remark. With respect to
our methodology, the mathematical results obtained by Pirah\~{a},
Cantor, and those presented in this paper do not contradict to
each other. \textit{They all are correct with respect to
mathematical languages used to express them.} This relativity is
very important and it has been emphasized in Postulate~2. For
instance, the result of Pirah\~{a} 1+2=`many' is correct in their
language in the same way as the result 1+2=3 is correct in the
modern mathematical languages. Analogously, the result
(\ref{4.4.1}) is correct in Cantor's language and the more
powerful language developed in this paper allows us to obtain a
more precise result (\ref{4.4.2}) that is correct in the new
language.

The choice of the mathematical language  depends on the practical
problem that are to be   solved and on the accuracy required for
such a solution. Again,  the result of Pirah\~{a} `many'+1=`many'
is correct. If one is satisfied with its accuracy, the answer
`many' can be used (and is used by Pirah\~{a}) in practice.
However, if one needs a more precise result, it is necessary to
introduce a more powerful mathematical language (a numeral system
in this case) allowing one to express the required answer in a
more accurate way.

\section{A brief conclusion}
\label{s7}

In this paper, a new computational methodology has been
introduced. It allows us to express, by a finite number of
symbols, not only finite numbers but infinite and infinitesimals,
too. All of them can be viewed as particular instances of a
general framework used to express numbers.

It  has been emphasized that the philosophical triad --
researcher, object of investigation, and tools used to observe the
object -- existing in such natural sciences as Physics and
Chemistry, exists in Mathematics, too. In natural sciences, the
instrument used to observe the object influences the results of
observations. The same happens in Mathematics where numeral
systems used to express numbers are among the instruments of
observations used by mathematicians. The usage of powerful numeral
systems gives the possibility to obtain more precise results in
Mathematics, in the same way as the usage of a good microscope
gives the possibility to obtain more precise results in Physics.

%\markboth{Bibliography}{Bibliography}
\bibliographystyle{amsplain}
\bibliography{XBib_numerical}
%\end{article}
\end{document}